\documentclass[12pt]{article}
\usepackage[T2A]{fontenc}
\usepackage[cp1251]{inputenc}
\usepackage[russian]{babel}
\usepackage{amsmath,amsthm,amsfonts}
\usepackage{amssymb,graphpap}

\newtheorem{lem}{Лемма}[section]
\newtheorem{theo}{Теорема}[section]
\newtheorem{cons}{Следствие}[theo]

\newtheorem{defn}{Определение}[section]

\newtheorem{exam}{Пример}[section]

\newtheorem{rem}{Примечание}[section]

\makeatletter
\newbox\para@box
\newcount\para@cnt
\def\ddet{%
\matrix@check\ddet \hbox\bgroup \everymath{}%
\setbox\para@box \hbox\bgroup $\env@matrix}
\def\endddet{%
\endmatrix
$\egroup \@tempdimb\dp\para@box \advance\@tempdimb by
-\ht\para@box \divide\@tempdimb by 3
\para@cnt\ht\para@box\relax
\advance\para@cnt by \dp\para@box\relax
\setbox\para@box\hbox{\raise\@tempdimb\box\para@box}%
\@tempdima=1mm \advance\para@cnt by \@tempdima
\@tempdima=1\unitlength\relax \@tempcnta \@tempdima\relax
\divide\para@cnt by \@tempcnta \divide\para@cnt by 6 \hbox{%
\unitlength=1\unitlength
$\hbox{\noindent\begin{picture}(\the\para@cnt,0)
\put(0,0){\line(1,3){\the\para@cnt}}
\put(0,0){\line(1,-3){\the\para@cnt}}
\end{picture}}\;%
\hbox{\noindent\box\para@box}\;%
\hbox{\noindent\begin{picture}(\the\para@cnt,0)
\put(\the\para@cnt,0){\line(-1,3){\the\para@cnt}}
\put(\the\para@cnt,0){\line(-1,-3){\the\para@cnt}}
\end{picture}}$%
}\egroup }

\advance\textwidth by 2cm \advance\textheight by 6cm \hoffset=-1cm
\begin{document}
\pagestyle{empty}

\def\brk#1{#1\discretionary{}{\hbox{$#1$}}{}}
\def\ge{\brk\geqslant} \let\geq=\ge
\def\le{\brk\leqslant} \let\leq=\le
\mathchardef\origplus=\mathcode`\+

\def\SP{\ignorespaces\hskip.5\parindent\relax}

\begin{center}

\Large \textbf{О рациональных приближениях алгебраических чисел высших порядков и некоторой параметризации  обобщенных уравнений Пелля.}
\vspace{6mm}

\large Заторский Р.А.
\end{center}
\vspace{6mm}
\footnotesize \hspace{6mm}Предложен новый алгебраический объект -- рекуррентные дроби, являющийся $n$-мерным обобщением непрерывных дробей, с помощью которого построен алгоритм рациональных приближений алгебраических иррациональностей. Для обобщенных уравнений Пелля построена некоторая параметризация.

\vspace{6mm}
 A new algebraic object is introduced - recurrent fractions, which is
an $n$-dimensional generalization of continued fractions.
It is used to describe an algorithm for rational approximations of
algebraic irrational numbers. Some parametrization for generalized
Pell's equations is constructed.
\vspace{6mm}

\section{Вступление}
\footnotesize1. \emph{Рациональные приближения алгебраических иррациональностей.}
Существует два вида алгоритмов связанных с алгебраическими иррациональностями высших порядков. К первому виду относятся алгоритмы построения по заданным алгебраическим иррациональностям $n$-мерных обобщений непрерывных дробей, а ко второму -- алгоритмы вычисления рациональных приближений к заданным $n$-мерным обобщениям непрерывных дробей. Обоим алгоритмам посвящено множество работ выдающихся аналитиков начиная с Эйлера.

  В основном, все усилия аналитиков были сконцентрированы в направлении обобщения непрерывных дробей и обобщения теоремы Лагранжа о рациональных приближениях квадратических иррациональностей укороченными периодическими цепными дробями. При этом самыми эффективными подходами оказались матричный, яркими представителями которого являются Эйлер\cite{Eu}, Пуанкаре\cite{Pu}, Якоби\cite{Ja}, Перрон\cite{Pe}, Брун\cite{Brun}, Бернштейн\cite{Bernst}, Пустыльников\cite{Pust}. Матричные алгоритмы, в общем случае, могут быть описаны рекуррентным равенством $$P_{k+1}=A_k\cdot P_k,$$  где $P_k$ -- $n$-мерный вектор, $A_k$ -- квадратная матрица с целыми элементами, некоторым образом, зависящая от вектора $P_k,$ причем $\text{det}A_k=\pm 1.$ Они просты, но не всегда обладают аналогом свойства периодичности цепных дробей для квадратических иррациональностей.

  Второй  подход  базируется на геометрии линейных однородных форм и целочисленных решеток, представителями которого являются Дирихле\cite{Di},  Ермит\cite{Er},  Клейн\cite{Kl}, Минковский\cite{Min}, Вороной\cite{Vor},  Скубенко\cite{Sk}, Арнольд\cite{Arn},  Брюно\cite{Br}. Эти алгоритмы дают хорошие приближения, но плохо программируются и обобщаются на высшие порядки.

  Рассмотрим отдельно два алгоритма имеющие непосредственное отношение к настоящей работе. Это алгоритмы Фюрстенау\cite{Furshtenau} и Пуанкаре\cite{Poankare}.

  Алгоритм Фюрстенау был предложен  в 1876 году и либо был незамечен либо забыт. В нем,  по заданным действительным числам $p$ и $q,$ строится последовательность равенств $$p=a_0+\frac{q_1}{p_1},\,p_1=a_1+\frac{q_2}{p_2},\,p_2=a_2+\frac{q_3}{p_3},\ldots,$$ $$q=b_0+\frac{c_1}{p_1},\,q_1=b_1+\frac{c_2}{p_2},\,q_2=b_2+\frac{c_3}{p_3},\ldots.$$
Далее, при помощи последовательных подстановок, строятся выражения для $p$ и $q:$

\begin{equation}\label{p}p=a_0+\frac{b_1+\frac{c_2}{a_2+\frac{b_3+
\frac{c_4}{a_4+\ldots}}{a_3+\frac{b_4+\ldots}{a_4+\ldots}}}}
{a_1+\frac{b_2+\frac{c_3}{a_3+\frac{b_4+\ldots}{a_4+\ldots}}}
{a_2+\frac{b_3+\frac{c_4}{a_4+\ldots}}{a_3+\frac{b_4+\ldots}{a_4+\ldots}}}},\end{equation}
\begin{equation}\label{q}q=b_0+\frac{c_1}{a_1+\frac{b_2+\frac{c_3}{a_3+
\frac{b_4+\ldots}{a_4+\ldots}}}{a_2+\frac{b_3+\frac{c_4}
{a_4+\ldots}}{a_3+\frac{b_4+\ldots}{a_4+\ldots}}}},\end{equation} которые Фюрстенау называет непрерывными дробями $2$-го класса. В случае $c_i=0,\,i=1,2,3,\ldots,$ они вырождаются в непрерывные дроби.

Суть алгоритма Пуанкаре, предложенного им в 1885 году, выражена  следующей теоремой:

\emph{Пусть последовательность $\{f_n\}_{n=0}^{\infty}$ является решением разностного уравнения $$f_{n+k}+\alpha_{n,1}f_{n+k-1}+\ldots+\alpha_{n,k}f_n=0,\,n=0,1,2,\ldots$$ с предельно постоянными коэффициентами, корни характеристического многочлена которого различны по модулю. Тогда либо $f_n=0$ при всех $n\geqslant n_0,$ либо существует предел $\underset{n\rightarrow\infty}{lim}\frac{f_{n+1}}{f_n},$ и этот предел равен одному из корней характеристического многочлена.}

В работе, для построения рациональных приближений алгебраических иррациональностей, использованы идеи алгоритмов Фюрстенау и Пуанкаре. Так как изображения (\ref{p}),(\ref{q}), которыми пользовался для своего обобщения Фюрстенау, не удобны в работе, то мы для их изображения  привлекаем аппарат параперманентов треугольных матриц \cite{ADM}.  При этом получаемые изображения не только близкие к изображениям непрерывных дробей, но и позволяют естественно ввести понятия порядка дроби и периодической дроби, а также с помощью параперманентов треугольных матриц, свойства которых хорошо изучены, показать, что периодические дроби высших порядков являются изображениями алгебраических иррациональностей высших порядков.

 Центральными теоремами в этой части работы являются теоремы \ref{theo.gran.rec.dr.=.korin} и   \ref{alg.rivn=rec.dr}. Эти теоремы устанавливают связь между максимальными по модулю действительными корнями алгебраических уравнений $n$-го порядка и рекуррентными дробями $n$-го порядка. Теорема \ref{forma-rivn} о связи некоторых алгебраических форм $n$-го порядка с алгебраическими уравнениями $n$-го порядка, взятая вместе с предыдущими двумя теоремами, в значительной степени решает вопрос о рациональных приближениях алгебраических иррациональностей.

2. \emph{Кольца целых алгебраических чисел полей $\mathbb{Q}(\sqrt[n]{m})$.}

Во второй части работы изучаются целые алгебраические числа вида $$s_0+s_1\sqrt[n]{m}+s_2\sqrt[n]{m^2}+\ldots+s_{n-1}\sqrt[n]{m^{n-1}}$$ и, в частности, единицы числовых полей $\mathbb{Q}(\sqrt[n]{m}).$ Отметим, что числовые поля $\mathbb{Q}(\sqrt[n]{m})$ при $n=2,3$ исследовались еще в работах  Вороного\cite{Voro}, Делоне и Фадеева\cite{Del}.

Описание структуры множества фундаментальных единиц в кольце целых чисел дает теорема Дирихле из которой следует, что по значению коэффициентов минимального многочлена при помощи конечного числа испытаний можно найти множество фундаментальных единиц.  Вороной построил алгоритм отыскания множества фундаментальных единиц кубических полей, но он сложный по числу операций и плохо обобщается на высшие порядки \cite{Buch}.

Известно\cite{Bernst}, что нахождение целых единиц поля $\mathbb{Q}(\sqrt[n]{m})$ связано с решением Диофантова  уравнения вида
$$\left|
  \begin{array}{cccccc}
    s_0 & ms_{n-1} & ms_{n-2} & \cdots & ms_2 & ms_1 \\
    s_1 & s_0 &ms_{n-1} & \cdots & ms_3 & ms_2 \\
    s_2 & s_1 & s_0 & \cdots & ms_4 & ms_3 \\
    \vdots & \cdots & \cdots & \cdots & \cdots & \vdots \\
    s_{n-2} & s_{n-3} & s_{n-4} & \cdots & s_0 & ms_{n-1} \\
    s_{n-1} & s_{n-2} & s_{n-3} & \cdots & s_1 & s_0 \\
  \end{array}
\right|=\pm 1.$$

Основным результатом второй части работы является некоторая параметризация этих уравнений при $$n=3,5,7,9,11.$$ При этом возникает числовой треугольник (см. примечание 5.1.). Дальнейшее изучение этого числового треугольника может пролить свет на общее решение приведенного выше диофантова уравнения.

Автор благодарен Манину Ю.И за консультации при написании этой работы.
\section{Предварительные сведения}

\normalsize Приведем некоторые сведения о парадетерминантах и параперманентах треугольных матриц, которые понадобятся нам в дальнейшем.

Пусть  $K$ ---  некоторое  поле.
\begin{defn}
\label{def.tr.matr} Треугольную таблицу \begin{equation}
\label{tr.matr} A= \left(
\begin{matrix}
a_{11}\\
a_{21} & a_{22} \\
\vdots & \vdots &\ddots\\
a_{n1} & a_{n2} & \cdots & a_{nn}
\end{matrix}
\right)_n
\end{equation} элементов поля
$K$ назовем \textbf{треугольной матрицей.}
\end{defn}
\begin{defn}
\label{defn.ugol} Каждому элементу $a_{ij}$
треугольной матрицы (\ref{tr.matr}) поставим в соответствие
треугольную матрицу с этим элементом в левом нижнем углу,
которую назовем \textbf{углом} заданной треугольной матрицы и
обозначим через $R_{ij}(A) $.
\end{defn}\index{ріг!треугольной матрицы}
Очевидно, что угол $R_{ij}(A) $ является треугольной матрицей
$(i - j + 1)$--го порядка. В угол $R_{ij}(A)$ входят только те
элементы $a_{rs}$ треугольной матрицы (\ref{tr.matr}), индексы
которых удовлетворяют соотношению $ j \leqslant s \leqslant r
\leqslant i$.

Ниже мы будем считать, что
\begin{equation*}
\label{ddet(R01)=1} \text{ddet}\, (R_{01}(A)) = \text{ddet}\,
(R_{n,n + 1}(A) ) = pper(R_{01}(A) ) = pper(R_{n,n + 1}(A) ) = 1 .
\end{equation*}

Так, в  треугольной матрице:
\begin{equation*}
\label{1.A5} A = \left( \begin{matrix}
 a_{11}&&&&\\
 a_{21}& a_{22}&&&\\
 a_{31}& a_{32}& a_{33}&&\\
 a_{41}& a_{42}& a_{43}& a_{44}&\\
 a_{51}& a_{52}& a_{53}& a_{54}& a_{55}\\
\end{matrix}  \right)
\end{equation*}
 угол $R_{42}(A) $  имеет вид: $$ R_{42}(A)  = \left(
\begin{matrix}
  a_{22}&&\\
  a_{32}& a_{33}&\\
  a_{42}& a_{43}& a_{44}\\
\end{matrix}  \right).
$$

\begin{defn}\emph{\cite{ADM}.}\label{theo.Gan}
 Если $A$ ---треугольная матрица (\ref{tr.matr}), то ее
 парадетерминантом и параперманентом назовем соответственно числа:
\begin{equation*}
\label{ddet2} \text{\em ddet}
(A)=\sum_{r=1}^{n}\sum_{p_{1}+\ldots+p_{r}=n}(-1)^{n-r}
\prod_{s=1}^{r}\{a_{p_{1}+\ldots+p_{s},p_{1}+\ldots+p_{s-1}+1}\},
\end{equation*}
$$
\text{\em pper(A)}=\sum_{r=1}^{n}\sum_{p_{1}+\ldots+p_{r}=n}
\prod_{s=1}^{r}\{a_{p_{1}+\ldots+p_{s},p_{1}+\ldots+p_{s-1}+1}\},
$$
где суммирование производится по множеству натуральных решений
уравнения $p_{1}+p_2+\ldots+p_{r}=n.$
\end{defn}

\begin{defn}
\label{def.tab} Прямоугольную таблицу элементов треугольной
матрицы (\ref{tr.matr}) назовем вписанной\index{вписанная
прямоугольная таблица} в эту матрицу, если одна ее вершина
совпадает с элементом $a_{n1}$, а противоположная к ней
--- с элементом  $a_{ii}$, $i \in \{1, \ldots ,n \}$. Эту таблицу
обозначим через $T(i)$.
\end{defn}

 Если  в определении \ref{def.tab}  $i = 1$, или $i = n$, то вписанная прямоугольная
 таблица вырождается соответственно в первый столбец или в $n$-тую
 строчку
  этой треугольной матрицы.

При нахождении значения парадетерминанта и параперманента
треугольных матриц удобно пользоваться \emph{алгебраическими
дополнениями}.
\begin{defn}
\label{def.alg} \textbf{Алгебраическими
дополнениями}\index{алгебраическое дополнение} $D_{ij}, P_{ij}
$ к факториальному произведению $\{a_{ij}\}$ ключевого элемента
$a_{ij}$ матрицы (\ref{tr.matr}) назовем соответственно числа
\begin{equation*}
\label{D(ij)=} D_{ij}=(-1)^{i+j}\cdot \text{\em ddet}\,
(R_{j-1,1})\cdot \text{\em ddet}\, (R_{n,i+1}),
\end{equation*}
\begin{equation*}
\label{P(ij)=} P_{ij}=pper(R_{j-1,1})\cdot pper(R_{n,i+1}),
\end{equation*}
где $R_{j-1,1}$ и $R_{n,i+1}$ --- углы треугольной матрицы
(\ref{tr.matr}).
\end{defn}

\begin{theo}
\label{theo.Rozklddet} \emph{\cite{ADM}. (Разложение
парафункции по элементам вписанной прямоугольной таблице).}
Пусть $A$
--- треугольная матрица (\ref{tr.matr}), а $T(i)$
--- некоторая вписанная в нее прямоугольная таблица элементов.
Тогда справедливо равенство:
\begin{eqnarray}\label{rozkl1}
&\text{\em ddet}\, (A)=\sum\limits_{s = 1}^i {\sum\limits_{r =
i}^n {\left\{ {a_{rs} } \right\}D_{rs} } },&\\\label{rozkl2}
&{pper}(A)=\sum\limits_{s = 1}^i {\sum\limits_{r = i}^n {\left\{
{a_{rs} } \right\}P_{rs} } },&
\end{eqnarray}
где $D_{rs} $ и $P_{rs}$ --- соответственно алгебраические
дополнения к факториальному произведению ключевого элемента
$a_{rs}$, который принадлежит таблице $T(i)$.
\end{theo}

\begin{cons}
\label{cons.vyrodtab}Если $i=1,$ то теорема
\ref{theo.Rozklddet} дает разложение парафункций по элементам
первого столбца и равенства (\ref{rozkl1}), (\ref{rozkl2})
примут вид:
\begin{equation*}
\label{rozk.i=1.d} \text{\em ddet}\, (A) = \sum\limits_{r = 1}^n
{\left\{ {a_{r1} } \right\}D_{r1} }  = \sum\limits_{r = 1}^n {( -
1)^{r + 1} \left\{ {a_{r1} } \right\} \cdot } \text{\em ddet}\,
(R_{n,r + 1} ) ,
\end{equation*}
\begin{equation*}
\label{rozk.i=1.p} pper(A) = \sum\limits_{r = 1}^n {\left\{
{a_{r1} } \right\}P_{r1} }  = \sum\limits_{r = 1}^n {\left\{
{a_{r1} } \right\} \cdot pper(R_{n,r + 1} )}.
\end{equation*}
Если же  $i=n$, то получим разложение парафункций по элементам
последней строчки:
\begin{equation*}
\label{rozk.i=n.d}\text{\em ddet}\, (A) = \sum\limits_{s = 1}^n
{\left\{ {a_{ns} } \right\}D_{ns} }  = \sum\limits_{s = 1}^n {( -
1)^{n + s} \left\{ {a_{ns} } \right\} \cdot }\text{\em ddet}\,
(R_{s - 1,1} ),
\end{equation*}
\begin{equation*}
\label{rozk.i=n.p}pper(A) = \sum\limits_{s = 1}^n {\left\{ {a_{ns}
} \right\}P_{ns} }  = \sum\limits_{s = 1}^n {\left\{ {a_{ns} }
\right\} \cdot pper(R_{s - 1,1} )}.
\end{equation*}
\end{cons}

Парадетерминанты и параперманенты треугольных матриц находят все больше применений в теории чисел и комбинаторном анализе. Более основательно с ними можно познакомится в \cite{ADM}.

\section{Рекуррентные дроби}

\begin{defn}\label{defn.rec.dr}
Пусть $a_{ij},\,1\leqslant j\leqslant i<\infty$--- некоторые
целые числа. Алгебраические объекты вида
$$\alpha=\frac{\left[\begin{array}{cccccccc}
  a_{11}  &  &  &  &  &  &  &    \\
  \frac{a_{22}}{a_{12}}  & a_{12} &  &  &  &  &  &    \\
  \vdots  & \ldots & \ddots &  &  &  &  &    \\
  \frac{a_{n-1,n-1}}{a_{n-2,n-1}} & \frac{a_{n-2,n-1}}{a_{n-3,n-1}} & \ldots & a_{1,n-1} &  &  &  &    \\
  \frac{a_{n,n}}{a_{n-1,n}} & \frac{a_{n-1,n}}{a_{n-2,n}} & \ldots & \frac{a_{2,n}}{a_{1,n}} & a_{1,n} &  &  &    \\
  0  & \frac{a_{n,n+1}}{a_{n-1,n+1}} & \ldots & \frac{a_{3,n+1}}{a_{2,n+1}} & \frac{a_{2,n+1}}{a_{1,n+1}} & a_{1,n+1} &  &   \\
  0  & 0 & \ldots & \frac{a_{4,n+2}}{a_{3,n+2}} & \frac{a_{3,n+2}}{a_{2,n+2}} & \frac{a_{2,n+2}}{a_{1,n+2}} & a_{1,n+2} &   \\
  \vdots & \ldots & \ldots & \ldots & \ldots & \ldots & \ldots & \ddots \\
\end{array}\right]_{\infty}}{\left[\begin{array}{ccccccc}
  a_{12} &  &  &  &  &  &    \\
 \vdots & \ddots &  &  &  &  &    \\
  \frac{a_{n-2,n-1}}{a_{n-3,n-1}} & \ldots & a_{1,n-1} &  &  &  &    \\
  \frac{a_{n-1,n}}{a_{n-2,n}} & \ldots & \frac{a_{2,n}}{a_{1,n}} & a_{1,n} &  &  &    \\
  \frac{a_{n,n+1}}{a_{n-1,n+1}} & \ldots & \frac{a_{3,n+1}}{a_{2,n+1}} & \frac{a_{2,n+1}}{a_{1,n+1}} & a_{1,n+1} &  &   \\
   0 & \ldots & \frac{a_{4,n+2}}{a_{3,n+2}} & \frac{a_{3,n+2}}{a_{2,n+2}} & \frac{a_{2,n+2}}{a_{1,n+2}} & a_{1,n+2} &   \\
 \vdots & \ldots & \ldots & \ldots & \ldots & \ldots & \ddots \\
\end{array}\right]_{\infty}}=$$\begin{equation}\label{rec.drib}=\left[\begin{array}{ccccccccc}
  a_{11} &  &  &  &  &  &  &  &    \\
  \frac{a_{22}}{a_{12}} & \vline & a_{12} &  &  &  &  &  &    \\
  \vdots & \vline & \ldots & \ddots &  &  &  &  &    \\
  \frac{a_{n-1,n-1}}{a_{n-2,n-1}} & \vline & \frac{a_{n-2,n-1}}{a_{n-3,n-1}} & \ldots & a_{1,n-1} &  &  &  &    \\
  \frac{a_{n,n}}{a_{n-1,n}} & \vline & \frac{a_{n-1,n}}{a_{n-2,n}} & \ldots & \frac{a_{2,n}}{a_{1,n}} & a_{1,n} &  &  &    \\
  0 & \vline & \frac{a_{n,n+1}}{a_{n-1,n+1}} & \ldots & \frac{a_{3,n+1}}{a_{2,n+1}} & \frac{a_{2,n+1}}{a_{1,n+1}} & a_{1,n+1} &  &   \\
  0 & \vline & 0 & \ldots & \frac{a_{4,n+2}}{a_{3,n+2}} & \frac{a_{3,n+2}}{a_{2,n+2}} & \frac{a_{2,n+2}}{a_{1,n+2}} & a_{1,n+2} &   \\
  \vdots & \vline & \ldots & \ldots & \ldots & \ldots & \ldots & \ldots & \ddots \\
\end{array}\right]_{\infty}\end{equation},

и
$$\frac{P_m}{Q_m}\!=\frac{\left[\begin{array}{ccccccccc}
  a_{11} &  &  &  &    &  &  &    \\
  \frac{a_{22}}{a_{12}} &  a_{12} &  &  &  &  &  &   \\
  \vdots & \ldots & \ddots &  &  &  &  &    \\
  \frac{a_{n,n}}{a_{n-1,n}} &  \frac{a_{n-1,n}}{a_{n-2,n}} & \ldots &  a_{1,n} &  &  &  &  \\
  0  & \frac{a_{n,n+1}}{a_{n-1,n+1}} & \ldots & \frac{a_{3,n+1}}{a_{2,n+1}} &  a_{1,n+1} &  &  & \\
  0  & 0 & \ldots & \frac{a_{4,n+2}}{a_{3,n+2}} & \frac{a_{3,n+2}}{a_{2,n+2}} &  a_{1,n+2} &  & \\
  \vdots  & \ldots & \ldots & \ldots &  \ldots & \ldots & \ddots& \\
  0 & 0 & \ldots & 0 & \frac{a_{n,m}}{a_{n-1,m}} & \frac{a_{n-1,m}}{a_{n-2,m}}   & \ldots & a_{1,m}\\
\end{array}\right]_{m}}{\left[\begin{array}{cccccccc}
  a_{12} &  &  &  &    &  &  &    \\
  \frac{a_{23}}{a_{13}} &  a_{13} &  &  &  &  &  &   \\
  \vdots & \ldots & \ddots &  &  &  &  &    \\
  \frac{a_{n-1,n}}{a_{n-2,n}} &  \frac{a_{n-2,n}}{a_{n-3,n}} & \ldots &  a_{1,n} &  &  &  &  \\
  \frac{a_{n,n+1}}{a_{n-1,n+1}}  & \frac{a_{n-1,n+1}}{a_{n-2,n+1}} &  \ldots &  \frac{a_{2,n+1}}{a_{1,n+1}} & a_{1,n+1} & & & \\
  0  & \frac{a_{n,n+2}}{a_{n-1,n+2}} & \ldots & \frac{a_{3,n+2}}{a_{2,n+2}} & \frac{a_{2,n+2}}{a_{1,n+2}} &  a_{1,n+2} &  & \\
  \vdots  & \ldots & \ldots & \ldots &  \ldots & \ldots & \ddots& \\
  0 & 0 & \ldots & 0 & \frac{a_{n,m}}{a_{n-1,m}} & \frac{a_{n-1,m}}{a_{n-2,m}}   & \ldots & a_{1,m}\\
\end{array}\right]_{m-1}}\!=$$\begin{equation}\label{rac.bk.rec.drib}=\left[\begin{array}{cccccccccc}
  a_{11} &  &  &  &  &  &  &  &    \\
  \frac{a_{22}}{a_{12}} & \vline & a_{12} &  &  &  &  &  &   \\
  \vdots & \vline & \ldots & \ddots &  &  &  &  &    \\
  \frac{a_{n,n}}{a_{n-1,n}} & \vline & \frac{a_{n-1,n}}{a_{n-2,n}} & \ldots &  a_{1,n} &  &  &  &  \\
  0 & \vline & \frac{a_{n,n+1}}{a_{n-1,n+1}} & \ldots & \frac{a_{3,n+1}}{a_{2,n+1}} &  a_{1,n+1} &  &  & \\
  0 & \vline & 0 & \ldots & \frac{a_{4,n+2}}{a_{3,n+2}} & \frac{a_{3,n+2}}{a_{2,n+2}} &  a_{1,n+2} &  & \\
  \vdots & \vline & \ldots & \ldots & \ldots &  \ldots & \ldots & \ddots& \\
  0 & \vline & 0 & \ldots & 0 & \frac{a_{n,m}}{a_{n-1,m}} & \frac{a_{n-1,m}}{a_{n-2,m}}   & \ldots & a_{1,m}\\
\end{array}\right]_{m},\end{equation}
 назовем соответственно
\textbf{рекуррентной дробью $n-$го порядка} и ее
\textbf{$m-$тым рациональным
укорочением}.\index{дробь!рекуррентная}\index{укорочение
рациональное}
\end{defn}

Разлагая параперманенты числителя и знаменателя рационального укорочения (\ref{rac.bk.rec.drib}) рекуррентной дроби (\ref{rec.drib})
 по элементам последней строчки, получим линейные рекуррентные
уравнения $n-$го порядка
\begin{equation*}\label{recQm=}
    P_m=a_{1m}P_{m-1}+a_{2m}P_{m-2}+\ldots+a_{nm}P_{m-n},\,\,m=1,2,\ldots,
\end{equation*}
\begin{equation*}\label{recPm=}
    Q_m=a_{1m}Q_{m-1}+a_{2m}Q_{m-2}+\ldots+a_{nm}Q_{m-n},\,\,m=1,2,\ldots,
\end{equation*}
где $$P_i=\begin{cases}1, \text{если}\,\, i=0,\\0,
\text{если}\,\, i<0,\end{cases}Q_i=\begin{cases}1,
\text{если}\,\, i=1-n,\\0, \text{если}\,\, 2-n\leqslant
i\leqslant 0,\end{cases}\,a_{n1}=1,$$ дающие эффективный
алгоритм вычисления значений рациональных укорочений
рекуррентных дробей $n-$го порядка.

Значение конечной границы $m-$тых рациональных укорочений
рекуррентной дроби $n-$го порядка, при $m\rightarrow\infty$
назовем значением соответственной рекуррентной дроби $n-$го
порядка.

Таким образом, рекуррентная дробь второго порядка при
$a_{1i}=q_i,\,a_{2i}=p_i$ имеет вид
$$\left[\begin{array}{cccccccc}
  q_1 &  &  &  &  &  &  &  \\
  \frac{p_2}{q_2} & \vline & q_2 &  &  &  &  &  \\
  0 & \vline & \frac{p_3}{q_3} & q_3 &  &  &  &  \\
  0 & \vline & 0 & \frac{p_4}{q_4} & q_4 &  &  &  \\
  \vdots & \vline & \ldots & \ldots & \ldots & \ddots &  &  \\
  0 & \vline & 0 & 0 & 0 & \frac{p_m}{q_m} &  &  \\
  \vdots & \vline & \ldots & \ldots & \ldots & \ldots & \ldots & \ddots \\
\end{array}\right]_{\infty}$$
и является иным изображением непрерывной дроби
$$q_1+\frac{p_2}{q_2+\frac{p_3}{q_3+\frac{p_4}{q_4+_{\ldots_{+_
{\frac{p_m}{q_m+_{\ldots}}}}}}}}.$$
\begin{defn}\label{def.rec.dr.per} Рекуррентную дробь $n-$го порядка
(\ref{rec.drib}) назовем $k-$пе\-ри\-одичес\-кой если ее
элементы удовлетворяют уравнениям
\begin{equation*}\label{rec.dr.per}
    a_{i,rk+j}=a_{i,j},\,i=1,2,\ldots,n;\,\,j=1,2,\ldots,k.
\end{equation*}
\end{defn}
\begin{defn}\label{defn.rivni.rec.dr}
Две рекуррентные дроби $n-$го порядка назовем равными, если
$m-$тые $(m=1,2,\ldots)$ рациональные укорочения этих дробей равны.
\end{defn}
\begin{defn}\label{rec.dr.zvych}
Рекуррентную дробь $n-$го порядка (\ref{rec.drib}) назовем
обыкновенной рекуррентной дробью, если выполнены равенства
$$a_{nn}=a_{n,n+1}=a_{n,n+2}=\ldots=1.$$
\end{defn}

\begin{lem}\label{lem.p_m/p_m-1}

i) Пусть $x_1,x_2,\ldots,x_n$ ---  различные корни алгебраического уравнения
\begin{equation}\label{algebr.rivn}
    x^n=a_1x^{n-1}+a_2x^{n-2}+\ldots+a_n,\,\,a_n\neq 0,
\end{equation}  причем модуль действительного корня $x_1$  превышает модули всех других корней этого уравнения, тогда
справедливо равенство
\begin{equation}\label{lim}\underset{m\rightarrow\infty}{\text{lim}}\frac{p_m(x_1,x_2,\ldots,x_n)}{p_{m-1}(x_1,x_2,\ldots,x_n)}=x_1,\end{equation} где
\begin{equation}\label{p_m(x_1...x_n)}p_m(x_1,x_2,\ldots,x_n)=\sum_{\alpha_1+\alpha_2+\ldots+\alpha_n=m}x_1^{\alpha_1}x_2^{\alpha_2}\cdot\ldots\cdot
x_n^{\alpha_n},\,\,m=1,2,\ldots -\end{equation}  полный однородный симметрический многочлен;

Справедливо и обратное утверждение. Если $x_1,x_2,\ldots,x_n$ ---  различные корни алгебраического уравнения (\ref{algebr.rivn}) и  справедливо равенство (\ref{lim}), то $x_1$ --  действительный корень наибольшего модуля алгебраического уравнения (\ref{algebr.rivn}).

ii) параперманент
 \begin{equation}\label{P_m}
 P_m=\left[\begin{array}{ccccccc}a_{1}&&&&&&\\
    \frac{a_{2}}{a_{1}}&a_{1}&&&&&\\
    \ldots&\ldots&\ddots&&&&\\
    \frac{a_{n}}{a_{n-1}}& \frac{a_{n-1}}{a_{n-2}}&\ldots&a_{1}&&&\\
     0& \frac{a_{n}}{a_{n-1}}&\ldots&\frac{a_{2}}
     {a_{1}}&a_{1}&&\\
     \ldots&\ldots&\ldots&\ldots&\ldots&\ddots&
     \\0&\ldots&0&\frac{a_{n}}{a_{n-1}}&\ldots&
     \frac{a_{2}}{a_{1}}&a_{1}\end{array}\right]_{m}, P_0=1
     \end{equation}
      и однородный симметрический многочлен $$p_m(x_1,x_2,\ldots,x_n),\,\, p_0(x_1,x_2,\ldots,x_n)=1$$ удовлетворяют одному и тому же линейному рекуррентному уравнению вида \begin{equation}\label{rec.ur.Pp}u_m=a_1u_{m-1}+a_2u_{m-2}+\ldots+a_nu_{m-n},\end{equation} т.е. их значения совпадают при любом натуральном $m.$
\end{lem}
\begin{proof}

i)

 Производящей функцией полного симметрического многочлена (\ref{p_m(x_1...x_n)}) является функция $\prod_{i=1}^n\frac{1}{1-x_iz}.$
 Разложим ее на простейшие дроби, получим $$\prod_{i=1}^n\frac{1}{1-x_iz}\!=\!\!\sum_{i=1}^n\frac{x_i^{n-1}}{(x_i-x_1)
 \!\cdot\ldots\cdot\!(x_i-x_{i-1})\!(x_i-x_{i+1})\!\cdot\ldots\cdot\!(x_i-x_n)\!(1-x_iz)}.$$
    Найдем коэффициент $p_m(x_1,x_2,\ldots,x_n)$ при $z^m$ в разложении последней функции в степенной ряд. Он равен $$\sum_{i=1}^n\frac{x_i^{m+n-1}}{(x_i-x_1)\cdot\ldots\cdot(x_i-x_{i-1})(x_i-x_{i+1})\cdot\ldots\cdot(x_i-x_n)}.$$ Следовательно $$\underset{m\rightarrow\infty}{\text{lim}}\frac{p_m(x_1,x_2,\ldots,x_n)}{p_{m-1}(x_1,x_2,\ldots,x_n)}=\underset{m\rightarrow\infty}
    {\text{lim}}\frac{\sum_{i=1}^n\frac{x_i^{m+n-1}}{(x_i-x_1)\cdot\ldots\cdot(x_i-x_{i-1})(x_i-x_{i+1})\cdot\ldots\cdot(x_i-x_n)}}
    {\sum_{i=1}^n\frac{x_i^{m+n-2}}{(x_i-x_1)\cdot\ldots\cdot(x_i-x_{i-1})(x_i-x_{i+1})\cdot\ldots\cdot(x_i-x_n)}}=$$
$$=x_1\cdot\underset{m\rightarrow\infty}{\text{lim}}\frac{1+\sum_{i=2}^n(\frac{x_i}{x_1})^{m+n-1}\frac{(x_1-x_2)(x_1-x_3)\cdot\ldots\cdot(x_1-x_n)}
{(x_i-x_1)\cdot\ldots\cdot(x_i-x_{i-1})(x_i-x_{i+1})\cdot\ldots\cdot(x_i-x_n)}}{1+\sum_{i=2}^n(\frac{x_i}{x_1})^{m+n-2}
\frac{(x_1-x_2)(x_1-x_3)\cdot\ldots\cdot(x_1-x_n)}{(x_i-x_1)\cdot\ldots\cdot(x_i-x_{i-1})(x_i-x_{i+1})\cdot\ldots\cdot(x_i-x_n)}}.$$
Таким образом, справедливость п. i) леммы \ref{lem.p_m/p_m-1} следует из последних равенств.

ii) То, что параперманент (\ref{P_m}) удовлетворяет линейному рекуррентному уравнению (\ref{rec.ur.Pp}) непосредственно видно из разложения этого параперманента по элементам первого столбца или последней строчки.

Докажем, что полный однородный симметрический многочлен также удовлетворяет этому рекуррентному уравнению.
Производящей функцией элементарных симметрических многочленов, построенных с помощью переменных $x_1,x_2,\ldots,x_n,$ является
$$\sigma(z)=\sum_{m=0}^{\infty}\sigma_mz^m=\prod_{i=0}^n(1+x_iz),\,\sigma_0=1,$$
а производящей функцией полных однородных многочленов, построенных с помощью тех же переменных,  как отмечалось выше, --- функция
$$p(z)=\sum_{m=0}^{\infty}p_mz^m=\prod_{i=0}^n(1-x_iz)^{-1},\,p_0=1.$$ Следовательно, справедливо равенство $$\sigma(z)p(-z)=1.$$
Таким образом, свертка $\sum_{i=0}^m(-1)^i\sigma_ip_{m-i}$ равна нулю, т.е. справедливо рекуррентное соотношение
$$p_m=\sigma_1p_{m-1}-\sigma_2p_{m-2}+\ldots+(-1)^{m-1}\sigma_m.$$ Учитывая теорему Виета для уравнения (\ref{algebr.rivn}) получим соотношения $\sigma_i=(-1)^{i-1}a_i$ и, следовательно, линейное рекуррентное уравнение
$$p_m=a_1p_{m-1}+a_2p_{m-2}+\ldots+a_np_{m-n}.$$
\end{proof}

\begin{theo}\label{theo.gran.rec.dr.=.korin} Пусть задано алгебраическое уравнение (\ref{algebr.rivn}) с попарно различными корнями.
Если для $m-$го рационального укорочения
$1-$пе\-ри\-оди\-ческой рекуррентной дроби $n-$го порядка
\begin{equation}\label{1per.rec.dr}
    \left[\begin{array}{ccccccccc}
      a_1 &  &  &  &  &  &  &  &  \\
      \frac{a_2}{a_1} & \vline & a_1 &  &  &  &  &  &  \\
      \vdots & \vline & \ldots & \ddots &  &  &  &  &  \\
      \frac{a_{n-1}}{a_{n-2}} & \vline & \frac{a_{n-2}}{a_{n-3}} & \ldots & a_1 &  &  &  &  \\
      \frac{a_n}{a_{n-1}} & \vline & \frac{a_{n-1}}{a_{n-2}} & \ldots & \frac{a_2}{a_1} & a_1 &  &  &  \\
      0 & \vline & \frac{a_n}{a_{n-1}} & \ldots & \frac{a_3}{a_2} & \frac{a_2}{a_1} & a_1 &  &  \\
      \vdots & \vline & \ldots & \ldots & \ldots & \ldots & \ldots & \ddots &  \\
      0 & \vline & 0 & \ldots & \frac{a_n}{a_{n-1}} & \frac{a_{n-1}}{a_{n-2}} & \frac{a_{n-2}}{a_{n-3}} & \ldots & a_1 \\
    \end{array}\right]_{\infty},
\end{equation} построенной при помощи коэффициентов этого уравнения,  существует конечный ненулевой действительный предел при $m\rightarrow\infty,$ т.е. $$\underset{m\rightarrow\infty}{\text{lim}}\frac{P_m}{Q_m}=x\neq 0,$$
то такая рекуррентная дробь $n-$того порядка является изображением
действительного корня алгебраического уравнения (\ref{algebr.rivn}),
  модуль которого больше всех модулей  других корней этого уравнения.
\end{theo}
\begin{proof} 1. Докажем сначала, что значение рекуррентной дроби (\ref{1per.rec.dr}) является корнем уравнения (\ref{algebr.rivn}). Разложим числитель $m-$го рационального укорочения
(\ref{1per.rec.dr}) рекуррентной дроби $n-$го порядка
 по элементам первого столбца:
$$P_m=a_1P_{m-1}+a_2P_{m-2}+\ldots+a_nP_{m-n}.$$ Поскольку $P_{m-1}=Q_m,$ то $$\frac{P_m}{Q_m}
=\frac{P_m}{P_{m-1}}=\frac{a_1P_{m-1}+a_2P_{m-2}+\ldots+a_nP_{m-n}}{P_{m-1}}=$$
$$=a_1+\frac{a_2}{\frac{P_{m-1}}{P_{m-2}}}+\frac{a_3}{\frac{P_{m-1}}{P_{m-3}}}+\ldots+\frac{a_n}{\frac{P_{m-1}}{P_{m-n}}}=$$
$$=a_1+\frac{a_2}{\frac{P_{m-1}}{P_{m-2}}}+\frac{a_3}{\frac{P_{m-1}}{P_{m-2}}
\cdot\frac{P_{m-2}}{P_{m-3}}}+\ldots+\frac{a_n}{\frac{P_{m-1}}{P_{m-2}}\cdot\frac{P_{m-2}}
{P_{m-3}}\cdot\ldots\cdot\frac{P_{m-n+1}}{P_{m-n}}}.$$ Так как по условию теоремы
$$\underset{m\rightarrow\infty}{\text{lim}}\frac{P_{m}}{O_{m}}=x\neq 0,$$
то получим уравнение
\begin{equation*}\label{x=a+b/x}x=a_1+\frac{a_2}{x}+\frac{a_3}{x^2}+\ldots+\frac{a_n}{x^{n-1}},\end{equation*}
или уравнение (\ref{algebr.rivn}).

2. Согласно  второй части леммы \ref{lem.p_m/p_m-1} числитель $m$-того рационального укорочения рекуррентной дроби совпадает с полным однородным симметрическим многочленом $m$-того порядка с $n$ переменными. Следовательно, справедливо равенство $$\underset{m\rightarrow\infty}{\text{lim}}\frac{P_m}{Q_m}=
\underset{m\rightarrow\infty}{\text{lim}}\frac{p_m(x_1,x_2,\ldots,x_n)}{p_{m-1}(x_1,x_2,\ldots,x_n)}=x,$$ где, согласно первой части леммы \ref{lem.p_m/p_m-1}, $x$ -- корень наибольшего модуля.
\end{proof}

Справедлива и обратная теорема:
\begin{theo}\label{alg.rivn=rec.dr}
Если алгебраическое уравнение (\ref{algebr.rivn}), все корни которого попарно различны, имеет действительный
корень $x_1$, модуль которого больше всех модулей других корней этого уравнения, то однопериодическая рекуррентная дробь
\begin{equation*}\label{1per.rec.dr1}
    \left[\begin{array}{ccccccccc}
      a_1 &  &  &  &  &  &  &  &  \\
      \frac{a_2}{a_1} & \vline & a_1 &  &  &  &  &  &  \\
      \vdots & \vline & \ldots & \ddots &  &  &  &  &  \\
      \frac{a_{n-1}}{a_{n-2}} & \vline & \frac{a_{n-2}}{a_{n-3}} & \ldots & a_1 &  &  &  &  \\
      \frac{a_n}{a_{n-1}} & \vline & \frac{a_{n-1}}{a_{n-2}} & \ldots & \frac{a_2}{a_1} & a_1 &  &  &  \\
      0 & \vline & \frac{a_n}{a_{n-1}} & \ldots & \frac{a_3}{a_2} & \frac{a_2}{a_1} & a_1 &  &  \\
      \vdots & \vline & \ldots & \ldots & \ldots & \ldots & \ldots & \ddots &  \\
      0 & \vline & 0 & \ldots & \frac{a_n}{a_{n-1}} & \frac{a_{n-1}}{a_{n-2}} & \frac{a_{n-2}}{a_{n-3}} & \ldots & a_1 \\
    \end{array}\right]_{\infty}
    \end{equation*} является изображением этого корня.
\end{theo}
\begin{proof} Запишем алгебраическое уравнение (\ref{algebr.rivn})
с ненулевым свободным членом в виде
\begin{equation*}\label{zah.rec.dr}
x=a_1+\frac{a_2+\frac{a_3+\frac{a_4+\cdots+\frac{a_n}{x}}{x}}{x}}{x}.
\end{equation*}
При помощи бесконечных вложений последнее уравнение можно записать в виде равенств
$$
x=a_1+\frac{a_2+\frac{a_3+\frac{a_4+\cdots+\frac{a_n}{a_1+\frac{a_2+\frac{a_3+\frac{a_4+\cdots+\frac{a_n}
{a_1+\ldots}}{a_1+\ldots}}{a_1+\ldots}}{a_1+\ldots}}}
{a_1+\frac{a_2+\frac{a_3+\frac{a_4+\cdots+\frac{a_n}{a_1+\ldots}}{a_1+\ldots}}{a_1+\ldots}}{a_1+\ldots}}}
{a_1+\frac{a_2+\frac{a_3+\frac{a_4+\cdots+\frac{a_n}{a_1+\ldots}}{a_1+\ldots}}{a_1+\ldots}}{a_1+\ldots}}}
{a_1+\frac{a_2+\frac{a_3+\frac{a_4+\cdots+\frac{a_n}{a_1+\ldots}}{a_1+\ldots}}{a_1+\ldots}}{a_1+\ldots}}=
$$\begin{equation}\label{x=zah.rec.dr}=a_1\vline+\frac{a_2\vline+\frac{a_3\vline+\frac{a_4\vline+\frac{a_5+\ldots}{a_1+\ldots}}
{a_1\vline+\frac{a_2+\ldots}{a_1+\ldots}}}{a_1\vline+\frac{a_2\vline+
\frac{a_3+\ldots}{a_1+\ldots}}{a_1\vline+\frac{a_2+\ldots}{a_1+\ldots}}}}{a_1\vline+\frac{a_2\vline+\frac{a_3\vline+\frac{a_4+\ldots}{a_1+\ldots}}
{a_1\vline+\frac{a_2+\ldots}
{a_1+\ldots}}}{a_1\vline+\frac{a_2\vline+\frac{a_3+\ldots}{a_1+\ldots}}{a_1\vline+\frac{a_2+\ldots}{a_1+\ldots}}}}.\end{equation}
правая сторона которых зависит только от коэффициентов уравнения\\
(\ref{zah.rec.dr}). Исследуем выражение правой части равенства (\ref{x=zah.rec.dr}).

 Первым
приближением  к значению этого выражения   является дробь, лежащая по левую сторону
первой вертикальной черточки, т.е.  дробь вида
$$\frac{a_1}{1}=\frac{u_1}{u_0}.$$ Другим приближением служит дробь, лежащая
слева второй вертикальной черточки
$$a_1+\frac{a_2}{a_1}=\frac{u_2}{u_1}=\frac{a_1u_1+a_2u_0}{a_1u_0}.$$ Третье приближение имеет вид
$$a_1+\frac{a_2+\frac{a_3}{a_1}}{a_1+\frac{a_2}{a_1}}=\frac{u_3}{u_2}=\frac{a_1u_2+a_2u_1+a_3u_0}{a_1u_1+a_2u_0}$$
и т. д.

По индукции можно показать, что $m$-тое
приближение имеет вид
$$\frac{a_1u_{m-1}+a_2u_{m-2}+\ldots+a_mu_0}{a_1u_{m-2}+a_2u_{m-3}+\ldots+a_{m-1}u_0},$$ если $m\leqslant n$ и
вид
$$\frac{a_1u_{m-1}+a_2u_{m-2}+\ldots+a_mu_{m-n}}{a_1u_{m-2}+a_2u_{m-3}+\ldots+a_mu_{m-n-1}},$$
если $m>n.$

Таким образом, $m$-тое приближение к значению исследуемого выражения совпадает с $m$-тым приближением  рекуррентной дроби вида
\begin{equation*}\label{1per.rec.dr1}
   \left[\begin{array}{ccccccccc}
      a_1 &&  &  &  &  &  &  &  \\
      \frac{a_2}{a_1}  &\vline& a_1 &  &  &  &  &  &  \\
      \vdots &\vline& \ldots & \ddots &  &  &  &  &  \\
      \frac{a_{n-1}}{a_{n-2}} &\vline& \frac{a_{n-2}}{a_{n-3}} & \ldots & a_1 &  &  &  &  \\
      \frac{a_n}{a_{n-1}}  &\vline& \frac{a_{n-1}}{a_{n-2}} & \ldots & \frac{a_2}{a_1} & a_1 &  &  &  \\
      0 &\vline& \frac{a_n}{a_{n-1}} & \ldots & \frac{a_3}{a_2} & \frac{a_2}{a_1} & a_1 &  &  \\
      \vdots  &\vline& \ldots & \ldots & \ldots & \ldots & \ldots & \ddots &  \\
      0  &\vline& 0 & \ldots & \frac{a_n}{a_{n-1}} & \frac{a_{n-1}}{a_{n-2}} & \frac{a_{n-2}}{a_{n-3}} & \ldots & a_1 \\
    \end{array}\right]_{\infty}.
    \end{equation*}
Но, согласно  лемме \ref{lem.p_m/p_m-1} справедливо равенство

$$\underset{m\rightarrow\infty}{\text{lim}}\left[\arraycolsep=2pt\begin{array}{ccccccccc}
      a_1 &&  &  &  &  &  &  &  \\
      \frac{a_2}{a_1}  &\vline& a_1 &  &  &  &  &  &  \\
      \vdots &\vline& \ldots & \ddots &  &  &  &  &  \\
      \frac{a_{n-1}}{a_{n-2}} &\vline& \frac{a_{n-2}}{a_{n-3}} & \ldots & a_1 &  &  &  &  \\
      \frac{a_n}{a_{n-1}}  &\vline& \frac{a_{n-1}}{a_{n-2}} & \ldots & \frac{a_2}{a_1} & a_1 &  &  &  \\
      0 &\vline& \frac{a_n}{a_{n-1}} & \ldots & \frac{a_3}{a_2} & \frac{a_2}{a_1} & a_1 &  &  \\
      \vdots  &\vline& \ldots & \ldots & \ldots & \ldots & \ldots & \ddots &  \\
      0  &\vline& 0 & \ldots & \frac{a_n}{a_{n-1}} & \frac{a_{n-1}}{a_{n-2}} & \frac{a_{n-2}}{a_{n-3}} & \ldots & a_1 \\
    \end{array}\right]_{m}\!\!\!\!\!=\!\!\underset{m\rightarrow\infty}{\text{lim}}\frac{p_m(x_1,x_2,\ldots,x_n)}{p_{m-1}(x_1,x_2,\ldots,x_n)}=x_1.$$

\end{proof}

\begin{exam}\label{exam.1}
\emph{Алгебраическое уравнение $$x^7=448x^6+672x^5+560x^4+280x^3+84x^2+14x+1$$ имеет действительный корень $$x=64+32\cdot \sqrt[7]{129}+16\cdot \sqrt[7]{129^2}+8\cdot \sqrt[7]{129^3}+4\cdot \sqrt[7]{129^4}+2\cdot \sqrt[7]{129^5}+ \sqrt[7]{129^6}\approx $$$$\approx449.49777653359235287015302078.$$ Найдем несколько первых рациональных укорочений соответственной рекуррентной дроби}
$$\left[\begin{array}{ccccccccccc}
  448 &&  &  &  &  &  &  &  &  &    \\
  3/2 &\vline& 448 &  &  &  &  &  &  &  &   \\
  5/6&\vline & 3/2 & 448 &  &  &  &  &  &  &    \\
  1/2 &\vline& 5/6 & 3/2 & 448 &  &  &  &  &  &    \\
  3/10 &\vline& 1/2 & 5/6 & 3/2 & 448 &  &  &  &  &    \\
  1/6 &\vline& 3/10 & 1/2 & 5/6 & 3/2 & 448 &  &  &  &    \\
  1/14&\vline & 1/6 & 3/10 & 1/2 & 5/6 & 3/2 & 448 &  &  &    \\
  0 &\vline& 1/14 & 1/6 & 3/10 & 1/2 & 5/6 & 3/2 & 448 &  &    \\
  0 &\vline& 0 & 1/14 & 1/6 & 3/10 & 1/2 & 5/6 & 3/2 & 448 &   \\
  \vdots &\vline& \ldots & \ldots & \ldots & \ldots & \ldots & \ldots & \ldots & \ldots & \ddots \\
\end{array}\right]_{\infty}
   $$
$$m=1: 448;$$
$$m=2: \frac{899}{2}=449.5;$$
$$m=3: \frac{808197}{1798}\approx 449.497775;$$
$$m=4: \frac{242188503}{538798}\approx 449.497776532;$$
$$m=5: \frac{217726387201}{484377006}\approx 449.497776533595; $$
$$m=6: \frac{195735053879083}{435452774402}\approx 449.497776533592351;$$
$$m=7: \frac{1231754601116629931}{2740290754307162}\approx 449.49777653359235286$$
$$m=8: \frac{553670954436947106368}{1231754601116629931}\approx 449.497776533592352870158;$$
$$m=9: \frac{1666566046544461900687}{3707617998461699373}\approx 449.4977765335923528701530208.$$
\emph{Таким образом, уже девятое рациональное укорочение рекуррентной дроби дает 24 верных десятичных знаков после точки.}
\end{exam}

\section{Алгебраические формы $n-$го порядка}
\begin{defn}\label{defn.forma}
\textbf{Алгебраической  $(n,m)$-формой} (далее
кратко\index{$(n,m)$-форма} \textbf{$(n,m)$-формой})  назовем
действительное
 число
 \begin{equation}\label{forma}
 x= F(n,m)=s_0+s_1\sqrt[n]{m}+\ldots+
 s_{n-1}\sqrt[n]{m^{n-1}},\,n\in\mathbb{N},\, s_i,m\in\mathbb{Q},
 \end{equation}
  или $n$-мерный вектор
  \begin{equation}\label{forma1}
  x=(s_0,s_1,\ldots,s_{n-1}).
  \end{equation}
\end{defn}

Очевидно, что множество $(n,m)$-форм с обычными операциями сложения и умножения образует поле.

\textbf{1. Изоморфизм $(n,m)$-форм с некоторыми классами
матриц}\index{изоморфизм!$(n,m)$-форм}

 $(n,m)$-форме (\ref{forma}) поставим в соответствие
циркулянт $n$-го порядка

$$X=$$\begin{equation}\label{cyrculant}\left(\arraycolsep=2pt
         \begin{array}{cccccc}
           s_0 & s_{n-1}\sqrt[n]{m^{n-1}} & s_{n-2}\sqrt[n]{m^{n-2}} & \cdots & s_2\sqrt[n]{m^2} & s_1\sqrt[n]{m} \\
           s_1\sqrt[n]{m} & s_0 & s_{n-1}\sqrt[n]{m^{n-1}} & \cdots & s_3\sqrt[n]{m^3} & s_2\sqrt[n]{m^2} \\
           s_2\sqrt[n]{m^2} & s_1\sqrt[n]{m} & s_0 & \cdots &  s_4\sqrt[n]{m^4} &s_3\sqrt[n]{m^3} \\
           \vdots & \cdots & \cdots & \cdots & \cdots & \vdots \\
           s_{n-2}\sqrt[n]{m^{n-2}} & s_{n-3}\sqrt[n]{m^{n-3}} & s_{n-4}\sqrt[n]{m^{n-4}} & \cdots & s_0 & s_{n-1}\sqrt[n]{m^{n-1}} \\
           s_{n-1}\sqrt[n]{m^{n-1}} & s_{n-2}\sqrt[n]{m^{n-2}} & s_{n-3}\sqrt[n]{m^{n-3}} & \cdots & s_1\sqrt[n]{m} & s_0 \\
         \end{array}
       \right),
\end{equation}  а $(n,m)$-форме (\ref{forma1})--- матрицу вида
\begin{equation}\label{cyrculant1}
X=\left(
  \begin{array}{cccccc}
    s_0 & ms_{n-1} & ms_{n-2} & \cdots & ms_2 & ms_1 \\
    s_1 & s_0 &ms_{n-1} & \cdots & ms_3 & ms_2 \\
    s_2 & s_1 & s_0 & \cdots & ms_4 & ms_3 \\
    \vdots & \cdots & \cdots & \cdots & \cdots & \vdots \\
    s_{n-2} & s_{n-3} & s_{n-4} & \cdots & s_0 & ms_{n-1} \\
    s_{n-1} & s_{n-2} & s_{n-3} & \cdots & s_1 & s_0 \\
  \end{array}
\right).
\end{equation}
Кождая из последних двух матриц однозначно задана своим первым
столбцом.

Поскольку произведение $(n,m)$-форм
\begin{equation}\label{x'}x'=s_0'+s_1'\sqrt[n]{m}+\ldots+s_{n-1}'\sqrt[n]{m^{n-1}},\end{equation}
\begin{equation}\label{x''}x''=s_0''+s_1''\sqrt[n]{m}+\ldots+s_{n-1}''\sqrt[n]{m^{n-1}}\end{equation}
является  $(n,m)$-формой вида
$$x=s_0+s_1\sqrt[n]{m}+\ldots+s_{n-1}\sqrt[n]{m^{n-1}},$$
где
\begin{equation}\label{dob.form}s_i=\sum_{j=0}^is_j's_{i-j}''+m\sum_{j=i+1}^{n-1}s_j's_{i-j}'',\,i=0,1,\ldots,n-1\end{equation}
и соответствует  первому столбцу произведения матриц $X'$ и $X'',$ которые соответствуют $(n,m)$-формам $(\ref{x'}), (\ref{x''}),$
то справедлива
\begin{theo}\label{izomorf} Множества $(n,m)$-форм $(\ref{forma}),$ $(\ref{forma1})$  изоморфны соответственно\index{изоморфизм}
множествам матриц $(\ref{cyrculant}),$ $(\ref{cyrculant1}).$
\end{theo}
Таким образом, $k$-той степени $(n,m)$-формы (\ref{x'})
 соответствует $k$-тая степень  матрицы
 $$X'=$$$$\left(\arraycolsep=2pt
         \begin{array}{cccccc}
           s_0' & s'_{n-1}\sqrt[n]{m^{n-1}} & s_{n-2}'\sqrt[n]{m^{n-2}} & \cdots & s_2'\sqrt[n]{m^2} & s_1'\sqrt[n]{m} \\
           s_1'\sqrt[n]{m} & s_0' & s_{n-1}'\sqrt[n]{m^{n-1}} & \cdots & s_3'\sqrt[n]{m^3} & s_2'\sqrt[n]{m^2} \\
           s_2'\sqrt[n]{m^2} & s_1'\sqrt[n]{m} & s_0' & \cdots &  s_4'\sqrt[n]{m^4} &s_3'\sqrt[n]{m^3} \\
           \vdots & \cdots & \cdots & \cdots & \cdots & \vdots \\
           s_{n-2}'\sqrt[n]{m^{n-2}} & s_{n-3}'\sqrt[n]{m^{n-3}} & s_{n-4}'\sqrt[n]{m^{n-4}} & \cdots & s_0' & s_{n-1}'\sqrt[n]{m^{n-1}} \\
           s_{n-1}'\sqrt[n]{m^{n-1}} & s_{n-2}'\sqrt[n]{m^{n-2}} & s_{n-3}'\sqrt[n]{m^{n-3}} & \cdots & s_1'\sqrt[n]{m} & s_0' \\
         \end{array}
       \right)$$ или матрицы
      $$ X'=\left(
  \begin{array}{cccccc}
    s_0' & ms_{n-1} & ms_{n-2}' & \cdots & ms_2' & ms_1' \\
    s_1' & s_0' &ms_{n-1}' & \cdots & ms_3' & ms_2' \\
    s_2' & s_1' & s_0' & \cdots & ms_4' & ms_3' \\
    \vdots & \cdots & \cdots & \cdots & \cdots & \vdots \\
    s_{n-2}' & s_{n-3}' & s_{n-4}' & \cdots & s_0' & ms_{n-1}' \\
    s_{n-1}' & s_{n-2}' & s_{n-3}' & \cdots & s_1' & s_0' \\
  \end{array}
\right).$$

Очевидно также, что если последние две матрицы умножить
соответственно на матрицы-столбцы
$$X''=\left(
                                                            \begin{array}{c}
                                                              s_0'' \\
                                                              s_1''\sqrt[n]{m} \\
                                                              s_2''\sqrt[n]{m^2} \\
                                                              \vdots \\
                                                              s_{n-2}''\sqrt[n]{m^{n-2}} \\
                                                              s_{n-1}''\sqrt[n]{m^{n-1}} \\
                                                            \end{array}
                                                          \right),\,X''=\left(
                                                                          \begin{array}{c}
                                                                            s_0'' \\
                                                                            s_1'' \\
                                                                            s_2'' \\
                                                                            \vdots \\
                                                                            s_{n-2}'' \\
                                                                            s_{n-1}'' \\
                                                                          \end{array}
                                                                        \right),
$$
то получим соответственно матрицы-столбцы\index{матрица!столбец}
$$X=\left(
                                                            \begin{array}{c}
                                                              s_0 \\
                                                              s_1\sqrt[n]{m} \\
                                                              s_2\sqrt[n]{m^2} \\
                                                              \vdots \\
                                                              s_{n-2}\sqrt[n]{m^{n-2}} \\
                                                              s_{n-1}\sqrt[n]{m^{n-1}} \\
                                                            \end{array}
                                                          \right),\,X=\left(
                                                                          \begin{array}{c}
                                                                            s_0 \\
                                                                            s_1 \\
                                                                            s_2 \\
                                                                            \vdots \\
                                                                            s_{n-2} \\
                                                                            s_{n-1} \\
                                                                          \end{array}
                                                                        \right),
$$ где $s_i$ задаются равенствами (\ref{dob.form}).

 Для любой $(n,m)$-формы\index{$(n,m)$-форма}
$$x=s_0+s_1\sqrt[n]{m}+\ldots+
 s_{n-1}\sqrt[n]{m^{n-1}}$$ можно найти лишь одну $(n,m)$-форму $$\overline{x}=\overline{s_0}+\overline{s_1}\sqrt[n]{m}+\ldots+
 \overline{s_{n-1}}\sqrt[n]{m^{n-1}}$$ такую, что их произведение
 $x\overline{x}$ является некоторым действительным числом. При этом
 $(n,m)$-форму $\overline{x}$ называют сопраженной
 $(n,m)$-формой к $(n,m)$-форме $x=F(n,m),$ а их произведение ---
 нормой\index{норма!$(n,m)$-формы}
 последей и обозначают через $|F(n,m)|$.

 Пусть $X$ и $\overline{X}$ --- матрицы соответственные
 $(n,m)$-форме $$x=s_0+s_1\sqrt[n]{m}+\ldots+s_{n-1}\sqrt[n]{m^{n-1}}$$ и сопряженной $(n,m)$-форме $\overline{x}.$
 Тогда $$X\cdot\overline{X}=|F(n,m)\cdot|E,$$ где $E$ --- единичная матрица.
 причем, норма $(n,m)$-формы $x$ равна детерминанту
 матрицы\index{детерминант}
 $X,$ а матрица соответственная сопряженной $(n,m)$-форме
 $\overline{x}$ является обратной матрицей к матрице $X,$ умноженной
 на детерминант матрицы $X.$

Таким образом, $n$-мерным
 обобщением уравнения Пелля\index{уравнение!Пелля}
 $$\left|
                                             \begin{array}{cc}
                                               s_0 & ms_1 \\
                                               s_1 & s_0 \\
                                             \end{array}
                                           \right|=s_0^2-ms_1^2=\pm
                                           1
$$ является уравнение
$$
\left|
  \begin{array}{cccccc}
    s_0 & ms_{n-1} & ms_{n-2} & \cdots & ms_2 & ms_1 \\
    s_1 & s_0 &ms_{n-1} & \cdots & ms_3 & ms_2 \\
    s_2 & s_1 & s_0 & \cdots & ms_4 & ms_3 \\
    \vdots & \cdots & \cdots & \cdots & \cdots & \vdots \\
    s_{n-2} & s_{n-3} & s_{n-4} & \cdots & s_0 & ms_{n-1} \\
    s_{n-1} & s_{n-2} & s_{n-3} & \cdots & s_1 & s_0 \\
  \end{array}
\right|=\pm 1.
$$

\section{Связь $(n,m)-$форм с алгебраическими уравнениями}
Найдем целые коэффициенты уравнения
\begin{equation}\label{zahalne.rivn.n}
 x^n=a_{n1}x^{n-1}+a_{n2}x^{n-2}+\ldots+a_{n,n-1}x^1+a_{n,n}
 \end{equation}
 корнем
которого является $(n,m)$-форма
$$
 x= F(n,m)=s_0+s_1\sqrt[n]{m}+\ldots+
 s_{n-1}\sqrt[n]{m^{n-1}},\,s_i\in\mathbb{Q},\,m\in \mathbb{N}.
$$

Пусть задана матрица
\begin{equation}\label{matrix.zahalna}X=\left(
      \begin{array}{cccc}
        a_{11} & a_{12} & \ldots & a_{1n} \\
        a_{21} & a_{22} & \ldots & a_{2n} \\
        \vdots & \cdots & \cdots & \vdots \\
        a_{n1} & a_{n2} & \cdots & a_{nn} \\
      \end{array}
    \right).
\end{equation}
Главный  минор  $r$-го порядка этой матрицы обозначим
через\index{минор} (см. \cite{Hantmacher}, стр. 13):
$$X\left(
                                   \begin{array}{cccc}
                                     i_1 & i_2 & \ldots & i_r \\
                                     i_1 & i_2 & \ldots & i_r \\
                                   \end{array}
                                 \right)
=\left|
                                                       \begin{array}{cccc}
                                                         a_{i_1,i_1} & a_{i_1,i_2} & \cdots & a_{i_1,i_r} \\
                                                         a_{i_2,i_1} & a_{i_2,i_2} & \cdots &  a_{i_2,i_r} \\
                                                         \vdots & \cdots & \cdots & \vdots \\
                                                          a_{i_r,i_1} & a_{i_r,i_2} & \cdots &  a_{i_r,i_r} \\
                                                       \end{array}
                                                     \right|,$$ где
$$i_1<i_2<\ldots<i_r.$$ Характеристическое уравнение\index{уравнение!характеристическое}
$$det(X-xE)=0,$$  матрицы  (\ref{matrix.zahalna}), как известно, имеет развернутый вид
\begin{equation*}\label{charakter.rivn}
x^n=\alpha_{n,1}x^{n-1}+\alpha_{n,2}x^{n-2}+\ldots+\alpha_{n,n-1}x^1+\alpha_{n,n},
\end{equation*}
где
\begin{equation}\label{koef.charakter.rivn}
\alpha_{n,j}=(-1)^{j-1}\sum_{1\leqslant i_1<i_2<\ldots<i_j\leqslant n}X\left(
                                                     \begin{array}{cccc}
                                                       i_1 & i_2 & \ldots & i_j \\
                                                       i_1 & i_2 & \ldots & i_j \\
                                                     \end{array}
                                                   \right).
\end{equation}
Согласно теореме
Гамильтона-Кели\index{теорема!Гамильтона-Кели}, каждая квадратная
матрица удовлетворяет своему характеристическому уравнению, поэтому
справедливо тождество
\begin{equation}\label{totojn.zahalna}
X^n=\alpha_{n,1}X^{n-1}+\alpha_{n,2}X^{n-2}+\ldots+\alpha_{n,n-1}X^1+\alpha_{n,n},
\end{equation} с коэффициентами (\ref{koef.charakter.rivn}),
где $X$ --- матрица (\ref{matrix.zahalna}).

Пусть матрица $X$ в  уравнении (\ref{totojn.zahalna}) задана
равенством (\ref{cyrculant1}), тогда коэффициенты $a_{n,j}$
уравнения (\ref{zahalne.rivn.n}), корнем которого является $(n,m)$-форма
(\ref{forma1}), можно найти пользуясь равенствами
(\ref{koef.charakter.rivn}). Таким образом, справедлива
\begin{theo}\label{forma-rivn}
Если $(n,m)$-форма $$x=s_0+s_1\sqrt[n]{m}+\ldots+
 s_{n-1}\sqrt[n]{m^{n-1}}$$ является корнем уравнения
 $$x^n=a_{n1}x^{n-1}+a_{n2}x^{n-2}+\ldots+a_{n,n-1}x^1+a_{n,n},$$
 то коэффициенты этого уравнения равны
$$a_{n,j}=(-1)^{j-1}\sum_{1\leqslant i_1<i_2<\ldots<i_j\leqslant n}X\left(
                                                     \begin{array}{cccc}
                                                       i_1 & i_2 & \ldots & i_j \\
                                                       i_1 & i_2 & \ldots & i_j \\
                                                     \end{array}
                                                   \right),$$ где $$X\left(
                                                     \begin{array}{cccc}
                                                       i_1 & i_2 & \ldots & i_j \\
                                                       i_1 & i_2 & \ldots & i_j \\
                                                     \end{array}
                                                   \right)$$
                                                   --- главные
                                                   миноры
                                                   матрицы $$X=\left(
  \begin{array}{cccccc}
    s_0 & ms_{n-1} & ms_{n-2} & \cdots & ms_2 & ms_1 \\
    s_1 & s_0 &ms_{n-1} & \cdots & ms_3 & ms_2 \\
    s_2 & s_1 & s_0 & \cdots & ms_4 & ms_3 \\
    \vdots & \cdots & \cdots & \cdots & \cdots & \vdots \\
    s_{n-2} & s_{n-3} & s_{n-4} & \cdots & s_0 & ms_{n-1} \\
    s_{n-1} & s_{n-2} & s_{n-3} & \cdots & s_1 & s_0 \\
  \end{array}
\right).$$
\end{theo}

Приведем формулы для нахождения коэффициентов
алгебраических уравнений, корнями которых
являются $(n,m)$-формы при $n=2,3,4,5.$
$$a_{21}=2s_0,\,a_{22}=-\left|
                               \begin{array}{cc}
                                 s_0 & ms_1 \\
                                 s_1 & s_0 \\
                               \end{array}
\right|;$$
$$a_{31}=3s_0,\,a_{32}=-3\left|
                               \begin{array}{cc}
                                 s_0 & ms_2 \\
                                 s_1 & s_0 \\
                               \end{array}
\right|,\,a_{33}=\left|
                   \begin{array}{ccc}
                     s_0 & ms_2 & ms_1 \\
                     s_1 & s_0 & ms_2 \\
                     s_2 & s_1 & s_o \\
                   \end{array}
                 \right|;
$$
$$a_{41}=4s_0,\,a_{42}=-4\left|
                               \begin{array}{cc}
                                 s_0 & ms_3 \\
                                 s_1 & s_0 \\
                               \end{array}
\right|-2\left|
                               \begin{array}{cc}
                                 s_0 & ms_2 \\
                                 s_2 & s_0 \\
                               \end{array}
\right|,$$$$a_{43}=4\left|
                   \begin{array}{ccc}
                     s_0 & ms_3 & ms_2 \\
                     s_1 & s_0 & ms_3 \\
                     s_2 & s_1 & s_0 \\
                   \end{array}
                 \right|,a_{44}=-\left|
                                    \begin{array}{cccc}
                                      s_0 & ms_3 & ms_2 & ms_1 \\
                                      s_1 & s_0 & ms_3 & ms_2 \\
                                      s_2 & s_1 & s_0 & ms_3 \\
                                      s_3 & s_2 & s_1 & s_0 \\
                                    \end{array}
                                  \right|.
                 $$
$$a_{51}=5s_0,\,a_{52}=-5\left|
                               \begin{array}{cc}
                                 s_0 & ms_4 \\
                                 s_1 & s_0 \\
                               \end{array}
\right|-5\left|
                               \begin{array}{cc}
                                 s_0 & ms_3 \\
                                 s_2 & s_0 \\
                               \end{array}
\right|,$$$$a_{53}=5\left|
                   \begin{array}{ccc}
                     s_0 & ms_4 & ms_3 \\
                     s_1 & s_0 & ms_4 \\
                     s_2 & s_1 & s_0 \\
                   \end{array}
                 \right|+5\left|
                   \begin{array}{ccc}
                     s_0 & ms_4 & ms_2 \\
                     s_1 & s_0 & ms_3 \\
                     s_3 & s_2 & s_0 \\
                   \end{array}
                 \right|,$$$$a_{54}\!=\!-5\left|
                                    \begin{array}{cccc}
                                      s_0\! & ms_4 & ms_3 & ms_2 \\
                                      s_1\! & s_0 & ms_4 & ms_3 \\
                                      s_2\! & s_1 & s_0 & ms_4 \\
                                      s_3\! & s_2 & s_1 & s_0 \\
                                    \end{array}
                                  \right|,a_{55}\!=\!\!\left|
                                                   \begin{array}{ccccc}
                                                     s_0 & ms_4 & ms_3 & ms_2 & ms_1 \\
                                                     s_1 & s_0 & ms_4 & ms_3 & ms_2 \\
                                                     s_2 & s_1 & s_0 & ms_4 & ms_3 \\
                                                     s_3 & s_2 & s_1 & s_0 & ms_4 \\
                                                     s_4 & s_3 & s_2 & s_1 & s_0 \\
                                                   \end{array}
                                                 \right|
                                  $$

Справедлива
\begin{theo}\label{theo.super+}
$(n,m^n+1)$-форма вида
\begin{equation*}\label{super.forma+}
m^{n-1}+m^{n-2}\sqrt[n]{m^n+ 1}+\ldots+m\sqrt[n]{(m^n+ 1)^{n-2}}+\sqrt[n]{(m^n+1)^{n-1}}
\end{equation*} является корнем алгебраического уравнения
\begin{equation*}\label{super.alg.rivn+}
x^n=\binom{n}{1}m^{n-1}x^{n-1}+\binom{n}{2}m^{n-2}x^{n-2}+\ldots+\binom{n}{n-1}mx+\binom{n}{n}
\end{equation*}
\end{theo}
\begin{proof}
Так как все главные миноры одинакового порядка
матрицы\index{минор!главный}
$$\left(
    \begin{array}{cccccc}
      m^{n-1} & m^n+ 1 & m(m^n+ 1) & \!\!\ldots\!\! & m^{n-3}(m^n+ 1) & m^{n-2}(m^n+ 1) \\
      m^{n-2} & m^{n-1} & m^n+ 1 & \!\!\ldots\!\! & m^{n-4}(m^n+ 1) & m^{n-3}(m^n+ 1) \\
      m^{n-3} & m^{n-2} & m^{n-1} & \!\!\ldots\!\! & m^{n-5}(m^n+ 1) & m^{n-4}(m^n+ 1) \\
      \cdots & \cdots & \cdots & \!\!\cdots\!\! & \cdots & \cdots \\
      m & m^2 & m^3 & \!\!\ldots\!\! & m^{n-1} & m^n+ 1 \\
      1 & m & m^2 & \!\!\ldots\!\! & m^{n-2} & m^{n-1} \\
    \end{array}
  \right)
$$ равны между
собой, то достаточно найти один из них.
Найдем главный минор $s$-того порядка матрицы
$$\left(
    \begin{array}{ccccc}
      m^{n-1} & m^n+ 1 & m(m^n+ 1) & \cdots & m^{s-2}(m^n+ 1) \\
      m^{n-2} & m^{n-1} & m^n+ 1 & \cdots & m^{s-3}(m^n+1) \\
      m^{n-3} & m^{n-2} & m^{n-1} & \cdots & m^{s-4}(m^n+ 1) \\
      \cdots & \cdots & \cdots & \cdots & \cdots \\
      m^{n-s} & m^{n-s+1} & m^{n-s+2} & \cdots & m^{n-1} \\
    \end{array}
  \right)
$$
Умножим первый столбец на $-m^r, r=1,2,\ldots,s-1$ и
сложим с $(r+1)$-вым столбцом, тогда получим
детерминант\index{детерминант} матрицы
$$\left(
    \begin{array}{ccccc}
      m^{n-1} & 1 &  m & \cdots & m^{s-2} \\
      m^{n-2} & 0 &1 & \cdots &  m^{s-3} \\
      m^{n-3} & 0 & 0 & \cdots &  m^{s-4} \\
      \cdots & \cdots & \cdots & \cdots & \cdots \\
      m^{n-s} & 0 & 0 & \cdots & 0 \\
    \end{array}
  \right).
$$
Разложим последний детерминант по элементам первого столбца,
получим $$(-1)^{s+1}m^{n-s} .$$

Таким образом, согласно теореме \ref{forma-rivn} коэффициент
$a_{n,s}$ равен $$(-1)^{s-1}(-1)^{s+1}m^{n-s} \binom{n}{s}=
m^{n-s} \binom{n}{s}.$$
\end{proof}
Приведем теорему о рациональном приближении $(n,m^n+1)$-форм.

\begin{theo}\label{rac.nabl.1}
$k$-тым рациональным приближением к $(n,m^n+1)$-форме является
выражение
$$\left[\begin{array}{ccccccc}
    a_{n1} &&  &  &  &  &  \\
    \frac{a_{n2}}{a_{n1}}&\vline & a_{n1} &  &  &  &  \\
    \vdots &\vline& \ldots & \ldots &  &  &  \\
    \frac{a_{nn}}{a_{n,n-1}} &\vline& \frac{a_{n,n-1}}{a_{n,n-2}} & \ldots & a_{n1} &  &  \\
    0 &\vline& \frac{a_{nn}}{a_{n,n-1}} & \ldots& \frac{a_{n2}}{a_{n1}} & a_{n1} &  \\
    \vdots &\vline& \ldots & \ldots & \ldots & \ldots & \ddots
  \end{array}\right]_k,
$$ где $a_{n1}=\binom{n}{1}m^{n-1},a_{ni}=\frac{\binom{n}{i}m^{n-i}}{\binom{n}{i-1}m^{n-i+1}}=\frac{n-i+1}{im},\,i=2,\ldots,n.$
\end{theo}
\begin{proof}
Эта теорема непосредственно следует из теоремы \ref{theo.super+},  теоремы \ref{alg.rivn=rec.dr} и теоремы Островского \cite{Prasolov} о том, что уравнение (\ref{zahalne.rivn.n}), все коэффициенты которого неотрицательны, причем наибольший общий делитель номеров положительных коэффициентов равен единице, имеет положительный корень, модуль которого больше модулей всех остальных корней этого уравнения.
\end{proof}

\begin{theo}\label{theo.super-}
$(n,m^n-1)$-форма вида
\begin{equation*}\label{super.forma-}
m^{n-1}+m^{n-2}\sqrt[n]{m^n-1}+\ldots+m\sqrt[n]{(m^n- 1)^{n-2}}+\sqrt[n]{(m^n-1)^{n-1}}
\end{equation*} является корнем алгебраического уравнения
$$
x^n=\binom{n}{1}m^{n-1}x^{n-1}-\binom{n}{2}m^{n-2}x^{n-2}+\ldots+(-1)^{n-2}\binom{n}{n-1}mx+$$
\begin{equation*}\label{super.alg.rivn-}+(-1)^{n-1}\binom{n}{n}
\end{equation*}
\end{theo}
\begin{proof} Аналогично доказательству теоремы \ref{theo.super+}
 с той лишь разницей, что минор $s$-того порядка имеет вид
$$\left|
    \begin{array}{ccccc}
      m^{n-1} & m^n- 1 & m(m^n- 1) & \cdots & m^{s-2}(m^n- 1) \\
      m^{n-2} & m^{n-1} & m^n-1 & \cdots & m^{s-3}(m^n-1) \\
      m^{n-3} & m^{n-2} & m^{n-1} & \cdots & m^{s-4}(m^n- 1) \\
      \cdots & \cdots & \cdots & \cdots & \cdots \\
      m^{n-s} & m^{n-s+1} & m^{n-s+2} & \cdots & m^{n-1} \\
    \end{array}
  \right|.
$$ Поэтому $$a_{ns}=(-1)^{s-1}m^{n-s}\binom{n}{s}.$$
\end{proof}
\section{Параметризация некоторых диофантовых уравнений}
Как уже отмечалось выше, при заданных
натуральных числах $n$ и
$m$ и обычных операциях суммы и произведения $(n,m)$-формы вида $s_0+s_1\sqrt[n]{m}+\ldots+
 s_{n-1}\sqrt[n]{m^{n-1}}$
образуют числовое поле. Единицы этих числовых полей, как отмечалось выше, являются решениями диофантова уравнения
\begin{equation}\label{uzah.Pella}
\left|
  \begin{array}{cccccc}
    s_0 & ms_{n-1} & ms_{n-2} & \cdots & ms_2 & ms_1 \\
    s_1 & s_0 &ms_{n-1} & \cdots & ms_3 & ms_2 \\
    s_2 & s_1 & s_0 & \cdots & ms_4 & ms_3 \\
    \vdots & \cdots & \cdots & \cdots & \cdots & \vdots \\
    s_{n-2} & s_{n-3} & s_{n-4} & \cdots & s_0 & ms_{n-1} \\
    s_{n-1} & s_{n-2} & s_{n-3} & \cdots & s_1 & s_0 \\
  \end{array}
\right|=\pm 1
\end{equation}
для $(n,m)\text{-формы}\,\,\,\, $$s_0+s_1\sqrt[n]{m}+\ldots+
 s_{n-1}\sqrt[n]{m^{n-1}}.$

Эти уравнения исключительно важны для изучения структуры групп единиц полей, образующих $(n,m)$-формами.

Приведем примеры диофантовых
уравнений\index{уравнение!диофантово}, аналогичных уравнению Пелля:

При $n=3$ диофантово уравнение (\ref{uzah.Pella}) в развернутом виде запишется так:
 \begin{equation}\label{pel.3}|F(3,m)|=s_0^3+s_1^3m+s_2^3m^2-3s_0s_1s_2m=1,\end{equation}
 его частичным решением являются
 \begin{equation*}\label{rozv.det.3-}\text{1)}\,\,\,
  m=\frac{k}{n}(p-3),s_0=(p-2)(p-1)-1,s_1=kn(p-2),s_2=n(p-1),
  \end{equation*}
\begin{equation}\label{rozv.det.3+}\text{2)}\,\,\, m=\frac{k}{n}(p+3),s_0=(p+2)(p+1)-1,s_1=kn(p+2),s_2=n(p+1),
\end{equation}
где $p=k^2n.$
 Сопряженные к приведенным выше решениям
  также являются  решениями этого уравнения.
Вот эти решения:
$$m=\frac{k}{n}(p-3),\overline{s_0}=1,\overline{s_1}=-nk,\overline{s_2}=n,$$
$$m=\frac{k}{n}(p+3),\overline{s_0}=1,\overline{s_1}=nk,\overline{s_2}=-n,$$
где $p=k^2n.$

Интересными оказываются также решения:
$$m=\frac{k}{n}\cdot\frac{p^4-6p^3+12p^2-9p+3}{(p-1)^3},$$$$s_0=(p-2)(p-1)-1,s_1=kn(p-2),s_2=n(p-1),$$
и
$$m=\frac{k}{n}\cdot\frac{p^4-6p^3+12p^2-9p+3}{(p-1)^3},s_0=1,s_1=-kn,s_2=n,$$
где $p=k^2n.$

Для нахождения целых единиц последовательных натуральных
чисел $m$ можно использовать формулы
$$\text{2)}\,\,\,m,\,\,\, s_0=(p-2)(p-1)-1,s_1=kn(p-2),s_2=n(p-1),$$
где $p=k^2n, n=\frac{3k}{k^3-m},$
$$\text{2)}\,\,\,m,\,\,\, s_0=(p-2)(p-1)-1,s_1=kn(p-2),s_2=n(p-1),$$
и $p=k^2n, n=\frac{3k}{k^3+m}.$

У \cite{Delone}, \cite{Wada H.} приведено
таблицы целых единиц для всех
натуральных значений не больших $70$ и $250$ соответственно.
Таблица Вада (H. Wada) была построена при помощи алгоритма Вороного. Однако, этот алгоритм плохо программируется, особенно для простых значений $m,$ и сложный.

\begin{exam}\label{gig}
\emph{При помощи параметризации (\ref{rozv.det.3+}) диофантова
уравнения (\ref{pel.3})  можно, например, получить его
численное решение для простого $900$-цифрового числа}
$$m=7239872283931086911936409478464844362964969622795111359$$
$$974934755582575468179641313855769259727641974609798138893255$$
$$525211182516975159833221705921278768192182515320441657882401$$
$$010319428973812885522960333217672876733752301462230992931472$$
$$909904767713450913240634917893148613016192814762213558096464$$
$$402885840589676123168302402549079960138122269110454764380482$$
$$620797212301828876999333286096636315349989543361123952286204$$
$$679568796728357647509618040687075931376908228422344437391718$$
$$598847301096960325148718064874625829071932050706560411935363$$
$$097129685477802808489613517468157406211390960108809283162275$$
$$086186714073920595021715659951203220403390183303468032708161$$
$$048148186272608066667799184997780960378585296561952602849633$$
$$261805881923550356499415702867592774593662323899734807289980$$
$$768201261654488304289541162853189580522455748182295602375207$$
$$298343703588386419599807672013076867243091847913861105093811$$
$$84003:$$

$$s_0=5241575068763353274345186852093150332602451917539901799$$
$$932282202278664005900121988404431292090932896640188251981680$$
$$326006813119732590132197563293690094647723171139535411516659$$
$$033189312863606609940147337891254811171282064557496488402243$$
$$617969658697073506591959974748435341221226735197582720381982$$
$$299223922948600684651564931177032848265679312201428827959001$$
$$983962562284328721926456192013117840486119368302870988948460$$
$$525587477841663879827970098875001133833949548499240589613091$$
$$535101060988495239283487157834796036977947636759176039961934$$
$$635983735619442324532910010967441297842323402115179325365224$$
$$245869281413952853346561463816495884323228588257539790111535$$
$$747953054994512900080325152302710423732719820889455301956641$$
$$861857282819629676959400601048120714983887659901278623613182$$
$$407879535162028668083576182989214435953539550841881883975282$$
$$254526281252669102943716630553173127420155080795347352070812$$
$$016381586784193433265072280590538909341829945733502452641625$$
$$931520640623887413956468162623695016290163104755307070450565$$
$$018386797664109367099985840845436742946338689084059438529370$$
$$980293569397369513919276599750601135570615318764192038306604$$
$$625048781009391683175362913228909806064167749424746695866250$$
$$167155583420998090725390393566830479390865308186951243348805$$
$$128955195908146359128484607204108672995078355407650230047302$$
$$525656852598329814864510572426666443716411444025574130235274$$
$$997855448580076587298245238361118200320906537084395631166706$$
$$552169049171921781816738304163317987099463732868320419686342$$
$$623468541037782483471350549375617711059947138423616726380877$$
$$798997774685802919400817546371507105036207154258163449863192$$
$$006072439921566891914553763235700570806637937053080647921850$$
$$659312149431138678443461371406707470566960647591054958666507$$
$$537418684478798007861256778237664233271272397518289447699995$$
$$52001,$$

$$s_1=5837382812453080192782510012686056895953207185324073163$$
$$326897981957933682949695513507983491029486920012188627504494$$
$$177893921657958220388784001370275856098037383085531762963522$$
$$971189156321683103476648884554121145707919627054822135024996$$
$$713598322792895735775812714508522934772101361352683861432720$$
$$073900114106912227058271263835488514173776210929713031728143$$
$$020658980959398029636120122487651614819173546990317107786517$$
$$178410215400856343858214716553426370423951480553096675712774$$
$$468958077699339832705487904142090168200408722770201141582664$$
$$560452449405642129440250357098685425446074566918060219791132$$
$$363007162746408115832954541627885798021635503947642925349684$$
$$282534014647216657434175460802228708011685170895150772325571$$
$$087615550292336931829309013572009092767009080325757523536591$$
$$510013992618754958068329997547475869531048142602447242118778$$
$$493806379124626745764202274260304140714663843554957895910898$$
$$698414186433127324323367579263127048764559073405306071771882$$
$$596826951621638034797468339758339899767973661216950764913598$$
$$416171091117071259067193567274604575185075976883699429429326$$
$$004590124578208042969022313006858521158939356428449850642416$$
$$658215888406381634954363972245970953773831113259190905145353$$
$$651910423387253629204803309274874145056261733316224414418819$$
$$882775119954109756433799795670377706668401680240225244798342$$
$$044866900658529842607199378616159767552297046673691480801264$$
$$171766248204246521165704944897310823728394813574848275027097$$
$$822090504131394457578987826385810634236567526497646123341539$$
$$87200,$$

$$s_2=6500915784301066706315708295484201254229616201498503033$$
$$225936513615050993420189616701616830178378071900836470490154$$
$$657732067457179957157657333794851256034894968136264330763326$$
$$483001404178215065556482015470095055770440336314987404707539$$
$$118687843607883638753785153045831999214567483283266963033324$$
$$468627959776406039326130905280422733774603240810130346764707$$
$$423355372082736401680347264707967374459504414866071427807368$$
$$389132735593909933021176828597417494084702604708765657062133$$
$$864315777701770554890487564019276485209166532788447229289032$$
$$521767451176009243507500541324892087509388127480599318372094$$
$$999301811585407932194506720014955341892574612453650759890353$$
$$844472637357540963976598058842306100267767040583926200220331$$
$$394772756556429720546343529771198353817635073590399931530726$$
$$700341972023621915539455287084102154147498148194464706290738$$
$$447968542450197421140538927732429020141632262305192523267917$$
$$425330324681548207507671518477644745204713871088804078414990$$
$$541919034200578212140280834001696317177767049206304425999422$$
$$159180067619206014434655679173528173324398774833309779609240$$
$$928488267568439044748816208777357041657525580606925980057317$$
$$088906832160714842155446850596898250005557236090621965750468$$
$$30440,$$
\emph{где} $n=k=\lfloor10^{300}\{10^{10}\pi\}\rfloor+374,$\emph{
а символом $\{\cdot\}$ обозначено дробную часть числа, т.е.}
$$k=n=89793238462643383279502884197169399375105820974944592307$$
$$816406286208998628034825342117067982148086513282306647093844$$
$$609550582231725359408128481117450284102701938521105559644622$$
$$948954930381964428810975665933446128475648233786783165271201$$
$$909145648566923460348610454326648213393607260249141273724587$$
$$0440$$
\emph{В этом примере $m$ --- простое $900$-цифровое число, а
$s_0,s_1,s_2$ --- числа, состоящие соответственно из $1800, 1500,
1300$ цифр.}
\end{exam}

При $n=5$  уравнение (\ref{uzah.Pella}) примет вид

$$|F(5,m)|=(s_{0}^5+s_{1}^5m+s_{2}^5m^2+s_{3}^5m^3+s_{4}^5m^4)-$$$$
-5(s_{0}^3s_{1}s_{4}m+s_{0}^3s_{2}s_{3}m+s_{1}^3s_{0}s_{2}m+s_{1}^3s_{3}s_{4}m^2+s_{2}^3s_{0}s_{4}m^2+
  +s_{2}^3s_{1}s_{3}m^2
$$$$+s_{3}^3s_{0}s_{1}m^2+s_{3}^3s_{2}s_{4}m^3
+s_{4}^3s_{0}s_{3}m^3+s_{4}^3s_{1}s_{2}m^3) +
5(s_{0}^2s_{1}^2s_{3}m+s_{0}^2s_{2}s_{4}^2m^2$$
$$+s_{1}s_{0}^2s_{2}^2m+s_{4}s_{0}^2s_{3}^2m^2+s_{0}s_{4}^2s_{1}^2m^2+s_{0}s_{2}^2s_{3}^2m^2
                  +s_{1}^2s_{4}s_{2}^2m^2+s_{2}s_{1}^2s_{3}^2m^2$$
$$
+s_{1}s_{3}^2s_{4}^2m^3+s_{4}^2s_{2}^2s_{3}m^3)-5s_{0}s_{1}s_{2}s_{3}s_{4}m^2=1.
$$

Его частичными решениями являются:
$$1)\, |F(5,m)|=1:m=\frac{k}{r}\cdot(rk^4-5),$$$$s_0=1,s_1=-rk^3,s_2=2rk^2,s_3=-2rk,s_4=r,$$
$$2)\, |F(5,m)|=1:m=\frac{k}{r}\cdot(rk^4-5),$$$$
s_0=1+30r^2k^8-25rk^4+r^4k^{16}-10r^3k^{12},$$$$
s_1=rk^3(r^3k^{12}-9r^2k^8+23rk^4-14),$$$$
s_2=rk^2(r^3k^{12}-8r^2k^8+17rk^4-7),$$$$
 s_3=kr(r^3k^{12}-7r^2k^8+12rk^4-3),$$$$ s_4=r(r^3k^{12}-6r^2k^8+8rk^4-1)$$

$$3)\, |F(5,m)|=1:m=\frac{k}{r}\cdot(rk^4+5),$$$$s_0=1,s_1=rk^3,s_2=-2rk^2,s_3=2rk,s_4=-r,$$
$$4)\, |F(5,m)|=1:m=\frac{k}{r}\cdot(rk^4+5),$$$$
s_0=1+30r^2k^8+25rk^4+r^4k^{16}+10r^3k^{12},$$$$
s_1=rk^3(r^3k^{12}+9r^2k^8+23rk^4+14),$$$$
s_2=rk^2(r^3k^{12}+8r^2k^8+17rk^4+7),$$$$
 s_3=kr(r^3k^{12}+7r^2k^8+12rk^4+3),$$$$ s_4=r(r^3k^{12}+6r^2k^8+8rk^4+1)$$

Запишем несколько взаимно сопряженных пар  решений диофантовых уравнений
(\ref{uzah.Pella}) при $n=7,9,11$:

$$1)\, |F(7,m)|=1:m=\frac{k}{r}\cdot(rk^6-7),$$$$s_0=1,s_1=-rk^5,s_2=3rk^4,s_3=-5rk^3,s_4=5rk^2,$$$$s_5=-3rk,s_6=r,$$
$$2)\, |F(7,m)|=1:m=\frac{k}{r}\cdot(rk^6-7),$$$$s_0=r^6k^{36}-21r^5k^{30}+161r^4k^{24}-539r^3k^{18}+721r^2k^{12}-245rk^6+1,$$$$
s_1=rk^5(r^5k^{30}-20r^4k^{24}+144r^3k^{18}-442r^2k^{12}+516rk^6-132),$$$$s_2=rk^4(r^5k^{30}-19r^4k^{24}+128r^3k^{18}-358r^2k^{12}+360rk^6-66),$$$$
s_3=rk^3(r^5k^{30}-18r^4k^{24}+113r^3k^{18}-286r^2k^{12}+244rk^6-30),$$$$s_4=rk^2(r^5k^{30}-17r^4k^{24}+99r^3k^{18}-225r^2k^{12}+160rk^6-12),$$$$
s_5=rk(r^5k^{30}-16r^4k^{24}+86r^3k^{18}-174r^2k^{12}+101rk^6-4),$$$$
s_6=r(r^5k^{30}-15r^4k^{24}+74r^3k^{18}-132r^2k^{12}+61rk^6-1)
$$

$$3)\, |F(7,m)|=1:m=\frac{k}{r}\cdot(rk^6+7),$$$$s_0=1,s_1=rk^5,s_2=-3rk^4,s_3=5rk^3,s_4=-5rk^2,$$$$s_5=3rk,s_6=-r,$$
$$4)\, |F(7,m)|=1:m=\frac{k}{r}\cdot(rk^6+7),$$$$s_0=r^6k^{36}+21r^5k^{30}+161r^4k^{24}+539r^3k^{18}+721r^2k^{12}+245rk^6+1,$$$$
s_1=rk^5(r^5k^{30}+20r^4k^{24}+144r^3k^{18}+442r^2k^{12}+516rk^6+132),$$$$s_2=rk^4(r^5k^{30}+19r^4k^{24}+128r^3k^{18}+358r^2k^{12}+360rk^6+66),$$$$
s_3=rk^3(r^5k^{30}+18r^4k^{24}+113r^3k^{18}+286r^2k^{12}+244rk^6+30),$$$$s_4=rk^2(r^5k^{30}+17r^4k^{24}+99r^3k^{18}+225r^2k^{12}+160rk^6+12),$$$$
s_5=rk(r^5k^{30}+16r^4k^{24}+86r^3k^{18}+174r^2k^{12}+101rk^6+4),$$$$
s_6=r(r^5k^{30}+15r^4k^{24}+74r^3k^{18}+132r^2k^{12}+61rk^6+1)
$$

$$1)\, |F(9,m)|=1:m=\frac{k}{r}\cdot(rk^8-3),$$$$s_0=1,s_1=-rk^7,s_2=rk^6,s_3=-rk^5,s_4=2rk^4,$$$$s_5=-2rk^3,s_6=rk^2,s_7=-rk,s_8=r,$$
$$2)\, |F(9,m)|=1:m=\frac{k}{r}\cdot(rk^8-3),$$$$s_0=9m^4k^{32}-54k^{24}m^3+97m^2k^{16}-48mk^8+1,$$
$$s_1=mk^7(9k^{24}m^3-51m^2k^{16}+84mk^8-35),$$$$s_2=mk^6(9k^{24}m^3-48m^2k^{16}+72mk^8-25),$$
$$s_3=k^5m(9k^{24}m^3-45m^2k^{16}+61mk^8-17),$$$$ s_4=mk^4(9k^{24}m^3-42m^2k^{16}+51mk^8-11),$$
$$s_5=k^3m(9k^{24}m^3-39m^2k^{16}+42mk^8-7), $$$$s_6=k^2m(9k^{24}m^3-36m^2k^{16}+34mk^8-4), $$
$$s_7=mk(9k^{24}m^3-33m^2k^{16}+27mk^8-2),$$$$s_8=m(9k^{24}m^3-30m^2k^{16}+21mk^8-1),$$

$$3)\, |F(9,m)|=1:m=\frac{k}{r}\cdot(rk^8+3),$$$$s_0=1,s_1=rk^7,s_2=-rk^6,s_3=rk^5,s_4=-2rk^4,$$$$s_5=2rk^3,s_6=-rk^2,s_7=rk,s_8=-r,$$
$$4)\, |F(9,m)|=1:m=\frac{k}{r}\cdot(rk^8+3),$$$$s_0=9m^4k^{32}+54k^{24}m^3+97m^2k^{16}+48mk^8+1,$$
$$s_1=mk^7(9k^{24}m^3+51m^2k^{16}+84mk^8+35),$$$$s_2=mk^6(9k^{24}m^3+48m^2k^{16}+72mk^8+25),$$
$$s_3=k^5m(9k^{24}m^3+45m^2k^{16}+61mk^8+17),$$$$ s_4=mk^4(9k^{24}m^3+42m^2k^{16}+51mk^8+11),$$
$$s_5=k^3m(9k^{24}m^3+39m^2k^{16}+42mk^8+7), $$$$s_6=k^2m(9k^{24}m^3+36m^2k^{16}+34mk^8+4), $$
$$s_7=mk(9k^{24}m^3+33m^2k^{16}+27mk^8+2),$$$$s_8=m(9k^{24}m^3+30m^2k^{16}+21mk^8+1),$$

Решениями диофантова уравнения
$$|F(11,m)|=$$$$\left|\arraycolsep=2pt
                                     \begin{array}{ccccccccccc}
                                       s_0 & ms_{10} & ms_9 & ms_8 & ms_7 & ms_6 & ms_5 & ms_4 & ms_3 & ms_2 & ms_1 \\
                                       s_1 & s_0 & ms_{10} & ms_9 & ms_8 & ms_7 & ms_6 & ms_5 & ms_4 & ms_3 & ms_2 \\
                                       s_2 & s_1 & s_0 & ms_{10} & ms_9 & ms_8 & ms_7 & ms_6 & ms_5 & ms_4 & ms_3 \\
                                       s_3 & s_2 & s_1 & s_0 & ms_{10} & ms_9 & ms_8 & ms_7 & ms_6 & ms_5 & ms_4 \\
                                       s_4 & s_3 & s_2 & s_1 & s_0 & ms_{10} & ms_9 & ms_8 & ms_7 & ms_6 & ms_5 \\
                                       s_5 & s_4 & s_3 & s_2 & s_1 & s_0 & ms_{10} & ms_9 & ms_8 & ms_7 & ms_6 \\
                                       s_6 & s_5 & s_4 & s_3 & s_2 & s_1 & s_0 & ms_{10} & ms_9 & ms_8 & ms_7 \\
                                       s_7 & s_6 & s_5 & s_4 & s_3 & s_2 & s_1 & s_0 & ms_{10} & ms_9 & ms_8 \\
                                       s_8 & s_7 & s_6 & s_5 & s_4 & s_3 & s_2 & s_1 & s_0 & ms_{10} & ms_9 \\
                                       s_9 & s_8 & s_7 & s_6 & s_5 & s_4 & s_3 & s_2 & s_1 & s_0 & ms_{10} \\
                                       s_{10} & s_9 & s_8 & s_7 & s_6 & s_5 & s_4 & s_3 & s_2 & s_1 & s_0 \\
                                     \end{array}
                                   \right|=1
$$ являются:

$$1)\, m=\frac{k}{r}\cdot(rk^{10}-11),$$$$s_0=1,s_1=-rk^9,s_2=5rk^8,s_3=-15rk^7,$$$$s_4=30rk^6,s_5=-42rk^5,s_6=42rk^4,s_7=-30rk^3,s_8=15rk^2,$$$$s_9=-5rk,s_{10}=r,$$
$$2)\, m=\frac{k}{r}\cdot(rk^{10}-11),$$$$s_0=r^{10}k^{100}-55r^9k^{90}+1265r^8k^{80}-15730r^7k^{70}+114037r^6k^{60}-$$$$-483637r^5k^{50}+1137015r^4k^{40}-1295910r^3k^{30}+527329r^2k^{20}-$$$$-32065rk^{10}+1,$$
$$s_1=rk^9(r^9k^{90}-54r^8k^{80}+1216r^7k^{70}-14749r^6k^{60}+103758r^5k^{50}-$$$$423776r^4k^{40}+947934r^3k^{30}-1005966r^2k^{20}+363493rk^{10}-16796),$$
$$s_2=rk^8(r^9k^{90}-53r^8k^{80}+1168r^7k^{70}-13811r^6k^{60}+94212r^5k^{50}-$$$$-370171r^4k^{40}+786526r^3k^{30}-774787r^2k^{20}+246779rk^{10}-8398),$$
$$s_3=rk^7(r^9k^{90}-52r^8k^{80}+1121r^7k^{70}-12915r^6k^{60}+85362r^5k^{50}-$$$$-322301r^4k^{40}+649346r^3k^{30}-591812r^2k^{20}+164814rk^{10}-3978),$$
$$s_4=rk^6(r^9k^{90}-51r^8k^{80}+1075r^7k^{70}-12060r^6k^{60}+77172r^5k^{50}-$$$$-279676r^4k^{40}+533294r^3k^{30}-448110r^2k^{20}+108134rk^{10}-1768),$$
$$s_5=rk^5(r^9k^{90}-50r^8k^{80}+1030r^7k^{70}-11245r^6k^{60}+69607r^5k^{50}-$$$$-241836r^4k^{40}+435590r^3k^{30}-336175r^2k^{20}+69589rk^{10}-728),$$
$$s_6=rk^4(r^9k^{90}-49r^8k^{80}+986r^7k^{70}-10469r^6k^{60}+62633r^5k^{50}-$$$$-208350r^4k^{40}+353750r^3k^{30}-249740r^2k^{20}+43849rk^{10}-273),$$
$$s_7=rk^3(r^9k^{90}-48r^8k^{80}+943r^7k^{70}-9731r^6k^{60}+56217r^5k^{50}-$$$$-178815r^4k^{40}+285563r^3k^{30}-183609r^2k^{20}+26998rk^{10}-91),$$
$$s_8=rk^2(r^9k^{90}-47r^8k^{80}+901r^7k^{70}-9030r^6k^{60}+50327r^5k^{50}-$$$$-152855r^4k^{40}+229069r^3k^{30}-133506r^2k^{20}+16204rk^{10}-26),$$
$$s_9=rk(r^9k^{90}-46r^8k^{80}+860r^7k^{70}-8365r^6k^{60}+44932r^5k^{50}-$$$$-130120r^4k^{40}+182538r^3k^{30}-95940r^2k^{20}+9454rk^{10}-6),$$
$$s_{10}=r(r^9k^{90}-45r^8k^{80}+820r^7k^{70}-7735r^6k^{60}+40002r^5k^{50}-$$$$-110285r^4k^{40}+144450r^3k^{30}-68085r^2k^{20}+5344rk^{10}-1),$$

$$3)\, m=\frac{k}{r}\cdot(rk^{10}+11),$$$$s_0=1,s_1=rk^9,s_2=-5rk^8,s_3=15rk^7,$$$$s_4=-30rk^6,s_5=42rk^5,s_6=-42rk^4,s_7=30rk^3,s_8=-15rk^2,$$$$s_9=5rk,s_{10}=-r,$$
$$4)\, m=\frac{k}{r}\cdot(rk^{10}+11),$$$$s_0=r^{10}r^{100}+55r^9k^{90}+1265r^8k^{80}+15730r^7k^{70}+114037r^6k^{60}+$$$$483637r^5k^{50}+1137015r^4k^{40}+1295910r^3k^{30}+527329r^2k^{20}+$$$$+32065rk^{10}+1,$$
$$s_1=rk^9(r^9k^{90}+54r^8k^{80}+1216r^7k^{70}+14749r^6k^{60}+103758r^5k^{50}+$$$$423776r^4k^{40}+947934r^3k^{30}+1005966r^2k^{20}+363493rk^{10}+16796),$$
$$s_2=rk^8(r^9k^{90}+53r^8k^{80}+1168r^7k^{70}+13811r^6k^{60}+94212r^5k^{50}+$$$$+370171r^4k^{40}+786526r^3k^{30}+774787r^2k^{20}+246779rk^{10}+8398),$$
$$s_3=rk^7(r^9k^{90}+52r^8k^{80}+1121r^7k^{70}+12915r^6k^{60}+85362r^5k^{50}+$$$$+322301r^4k^{40}+649346r^3k^{30}+591812r^2k^{20}+164814rk^{10}+3978),$$
$$s_4=rk^6(r^9k^{90}+51r^8k^{80}+1075r^7k^{70}+12060r^6k^{60}+77172r^5k^{50}+$$$$+279676r^4k^{40}+533294r^3k^{30}+448110r^2k^{20}+108134rk^{10}+1768),$$
$$s_5=rk^5(r^9k^{90}+50r^8k^{80}+1030r^7k^{70}+11245r^6k^{60}+69607r^5k^{50}+$$$$+241836r^4k^{40}+435590r^3k^{30}+336175r^2k^{20}+69589rk^{10}+728),$$
$$s_6=rk^4(r^9k^{90}+49r^8k^{80}+986r^7k^{70}+10469r^6k^{60}+62633r^5k^{50}+$$$$+208350r^4k^{40}+353750r^3k^{30}+249740r^2k^{20}+43849rk^{10}+273),$$
$$s_7=rk^3(r^9k^{90}+48r^8k^{80}+943r^7k^{70}+9731r^6k^{60}+56217r^5k^{50}+$$$$+178815r^4k^{40}+285563r^3k^{30}+183609r^2k^{20}+26998rk^{10}+91),$$
$$s_8=rk^2(r^9k^{90}+47r^8k^{80}+901r^7k^{70}+9030r^6k^{60}+50327r^5k^{50}+$$$$+152855r^4k^{40}+229069r^3k^{30}+133506r^2k^{20}+16204rk^{10}+26),$$
$$s_9=rk(r^9k^{90}+46r^8k^{80}+860r^7k^{70}+8365r^6k^{60}+44932r^5k^{50}+$$$$+130120r^4k^{40}+182538r^3k^{30}+95940r^2k^{20}+9454rk^{10}+6),$$
$$s_{10}=r(r^9k^{90}+45r^8k^{80}+820r^7k^{70}+7735r^6k^{60}+40002r^5k^{50}+$$$$+110285r^4k^{40}+144450r^3k^{30}+68085r^2k^{20}+5344rk^{10}+1)$$
и т. д.

\begin{rem}\label{odynyci.1}
Модули коэффициентов многочленов, которые являются решениями диофантовых
уравнений $|F(2n-1,m)|=1,$ очевидно связаны с
числовым\index{треугольник!числовой} треугольником
$$\label{tr.odynyc}\begin{array}{ccccccc}
        n=2: && 1 &  &  &  &  \\
        n=3: && 1 & 2 &  &  &  \\
        n=4: && 1 & 3 & 5 &  &  \\
        n=5: && 1 & 1 & 1 & 2 &  \\
        n=6: && 1 & 5 & 15 & 30 & 42
      \end{array},
$$ а элементы этого числового треугольника с факторизацией числа $2n-1.$
\end{rem}

\begin{rem}\label{odynyci.2}
Свободные члены многочленов, которые являются решениями 2) и 4) также
связаны с приведенным выше числовым треугольником.

Например, для модулей свободных членов
$s_6^0,s_7^0,s_8^0,s_9^0,s_{10}^0$ многочленов
$s_6,s_7,s_8,s_9,s_{10}$ решений 2) и 4)
диофантова уравнения $|F(11,m)|=1$
справедливы соотношения:
$$s_{10}^0=1,$$ $$s_9^0-s_{10}^0=5,$$
$$s_8^0-2s_9^0+s_{10}^0=15,$$ $$s_7^0-3s_8^0+3s_9^0-s_{10}^0=30,$$
$$s_6^0-4s_7^0+6s_8^0-4s_9^0+s_{10}^0=42.$$
\end{rem}

Приведем примеры решений уравнения
$$|F(n,m)|=\pm 1,$$ которые вытекают из теорем \ref{theo.super+},
\ref{theo.super-}.
\begin{theo}\label{F=1}
 Решениями диофантовых решений
$$|F(n,m^n-1)|=1, n=2,3,\ldots,$$
$$|F(2n-1,m^{2n-1}+1)|=1, n=2,3,\ldots,$$
$$|F(2n,m^{2n}+1)|=-1, n=1,2,\ldots$$
 являются соответственно:
$$s_0=m^{n-1},s_1=m^{n-2},\ldots,s_{n-2}=m,s_{n-1}=1,$$
$$s_0=m^{2n-2},s_1=m^{2n-3},\ldots,s_{2n-3}=m,s_{2n-2}=1,$$
$$s_0=m^{2n-1},s_1=m^{2n-2},\ldots,s_{2n-2}=m,s_{2n-1}=1.$$
\end{theo}
Легко находятся и сопряженные решения к приведенным в этой теореме
решениям. \vspace{5 cm}


\emph{Прикарпатский национальный университет им. В. Стефаника}

\emph{Украина}

romazz@rambler.ru
\end{document}